\documentclass[12pt,a4paper,reqno]{amsart}
\usepackage{amsmath}
\usepackage{amsfonts}
\usepackage{amssymb}

\numberwithin{equation}{section}

\theoremstyle{plain}

\newtheorem{theorem}[subsubsection]{Theorem}
\newtheorem{proposition}[subsubsection]{Proposition}
\newtheorem{lemma}[subsubsection]{Lemma}

\newtheorem{conjecture}[subsubsection]{Conjecture}

\newtheorem*{dickson}{Dickson's Conjecture}

\theoremstyle{definition}

\newtheorem{definition}[subsubsection]{Definition}

\newtheorem{question}[subsubsection]{Question}

\renewcommand{\leq}{\leqslant}
\renewcommand{\geq}{\geqslant}

\newsavebox{\proofbox}
\savebox{\proofbox}{\begin{picture}(7,7)%
  \put(0,0){\framebox(7,7){}}\end{picture}}

\newcommand{\md}[1]{\ensuremath{(\operatorname{mod}\, #1)}}
\newcommand{\mdsub}[1]{\ensuremath{(\mbox{\scriptsize mod}\, #1)}}
\newcommand{\mdlem}[1]{\ensuremath{(\mbox{\textup{mod}}\, #1)}}
\newcommand{\mdsublem}[1]{\ensuremath{(\mbox{\scriptsize \textup{mod}}\, #1)}}

\newcommand{\ssubsection}[1]{%
     \subsection[#1]{\sc #1}}

\newcommand\E{\mathbb{E}}
\newcommand\Z{\mathbb{Z}}
\newcommand\R{\mathbb{R}}

\newcommand\C{\mathbb{C}}
\newcommand\N{\mathbb{N}}

\newcommand\vol{\operatorname{vol}}

\newcommand\eps{\varepsilon}

\def\proof{\noindent\textit{Proof. }}

\def\remarks{\noindent\textit{Remarks. }}
\def\endproof{\hfill{\usebox{\proofbox}}}

\def\yildirim{Y{\i}ld{\i}r{\i}m }
\def\half{\textstyle\frac{1}{2}\displaystyle}

\begin{document}

\title{Three topics in additive prime number theory}

\author{Ben Green}
\address{Centre for Mathematical Sciences\\
Wilberforce Road\\
     Cambridge CB3 0WA\\
     England
}
\email{b.j.green@dpmms.cam.ac.uk}

\subjclass{}

\begin{abstract}
We discuss, in varying degrees of detail, three contemporary themes in prime number theory. Topic 1: the work of Goldston, Pintz and \yildirim on short gaps between primes. Topic 2: the work of Mauduit and Rivat, establishing that 50\% of the primes have odd digit sum in base 2. Topic 3: work of Tao and the author on linear equations in primes.
\end{abstract}

\maketitle

\textsc{Introduction.} These notes are to accompany two lectures I am scheduled to give at the \emph{Current Developments in Mathematics} conference at Harvard in November 2007. The title of those lectures is `A good new millennium for primes', but I have chosen a rather drier title for these notes for two reasons. Firstly, the title of the lectures was unashamedly stolen (albeit with permission) from Andrew Granville's entertaining article \cite{granville-article} of the same name. Secondly, and more seriously, I do not wish to claim that the topics chosen here represent a complete survey of developments in prime number theory since 2000 or even a selection of the most important ones. Indeed there are certainly omissions, such as the lack of any discussion of the polynomial-time primality test \cite{aks}, the failure to even mention the recent work on primes in orbits by Bourgain, Gamburd and Sarnak, and many others.

I propose to discuss the following three topics, in greatly varying degrees of depth. Suggestions for further reading will be provided. The three sections may be read independently although there are links between them.

1. Gaps between primes. Let $p_n$ be the $n$th prime number, thus $p_1 = 2$, $p_2 = 3$, and so on. The prime number theorem, conjectured by Gauss and proven by Hadamard and de la Vall\'ee Poussin over 110 years ago, tells us that $p_n$ is asymptotic to $n\log n$, or in other words that
\[ \lim_{n \rightarrow \infty} \frac{p_n}{n\log n} = 1.\]
This implies that the gap between the $n$th and $(n+1)$st primes, $p_{n+1} - p_n$, is about $\log n$ on average. About 2 years ago Goldston, Pintz and \yildirim proved the following remarkable result: for any $\epsilon > 0$, there are infinitely many $n$ such that $p_{n+1} - p_n < \epsilon \log n$. That is, infinitely often there are consecutive primes whose spacing is \emph{much} closer than the average.

2. Digits of primes. Written in binary, the first few primes are 
\[ 10, 11, 101, 111, 1011, 1101, 10001, 10011, 10111,\dots\]
There is no obvious pattern\footnote{Except, of course, that the last digit of primes except the first is always 1.}. Indeed, why would there be, since the definition of `prime' has nothing to do with digital expansions. Proving such a statement, or even formulating it correctly, is an entirely different matter. A couple of years ago, however, Mauduit and Rivat did manage to prove that the digit sum is odd 50\% of the time (and hence even 50\% of the time). They also obtained results in other bases. 

3. Patterns of primes. Additive questions concerning primes have a long history. It has been known for over 70 years that there are infinitely many 3-term arithmetic progressions of primes such as $3,5,7$ and $5,11,17$, and that every large odd number is the sum of three primes. Recently, in joint work with Tao, we have been able to study more complicated patterns of primes. In this section we provide a guide to this recent joint work.

Throughout these notes we will write
\[ \E_{x \in X} f(x) := \frac{1}{X}\sum_{x \in X} f(x),\] where $X$ is any finite set and $f : X \rightarrow \C$ is a function.

\section{Gaps between primes}
\label{sec1}
These notes were originally prepared for a series of lectures I gave at the Norwegian Mathematical Society's \emph{Ski og mathematikk}, which took place at Rondablikk in January 2006. It is a pleasure to thank Christian Skau for inviting me to that event. The argument of Goldston, Pintz and \yildirim was first described to me by K.~Soundararajan at the Highbury Vaults in Bristol. It is a pleasure to thank him, and to refer the interested reader to his lectures on the subject \cite{soundararajan}, which are superior to these in every respect.

\ssubsection{The result}
     
   In 2005 Goldston, Pintz and \yildirim created a sensation by announcing a proof that
   \[ \mbox{liminf}_{n \rightarrow \infty} \frac{p_{n+1} - p_n}{\log n} = 0,\]
   where $p_n$ denotes the $n$th prime number. According to the prime number theorem we have
   \[ p_n \sim n \log n,\]
   and therefore
   \[ \frac{p_{n+1} - p_n}{\log n} \] has average value $1$. The Goldston, Pintz and \yildirim result thus states that the distance between consecutive primes can be $\epsilon$ of the average spacing, for any $\epsilon$, and is thus certainly most spectacular. 
   
   Previous efforts at locating small gaps between primes focussed on proving successively smaller upper bounds for $C := \mbox{liminf}_{n \rightarrow \infty} \frac{p_{n+1} - p_n}{\log n}$. The following table describing the history of these improvements does not make the Goldston-Pintz-\yildirim result look any less striking:
   
   \begin{tabular}{|l|l|l|l|}
   \hline & & $C$ & \\ \hline
   Trivial from PNT & & 1 & \\
   Hardy-Littlewood \cite{HL-0} & 1926 & 2/3 & on GRH \\
  Rankin \cite{Rankin} & & 3/5 & on GRH \\
Erd\H{o}s \cite{Erdos} & 1940 & 1 - c & unconditionally\\
Ricci \cite{Ricci} & 1954 & 15/16 & \\ 
 Bombieri-Davenport \cite{BD} & 1965 &  $0.4665\dots$ & \\ 
Pilt'ai \cite{Pi} & 1972 &  $0.4571\dots$ & \\ Uchiyama \cite{Uc} & & $0.4542\dots$ & \\
 Huxley \cite{Hu1,Hu2,Hu3} & 1984 &  $0.4393\dots$ & \\
 Maier \cite{Ma} & 1989 & $0.2484\dots$ & \\
 Goldston-Pintz-\yildirim & 2005 & 0 & \\
 \hline
   \end{tabular} 
   \vspace{20pt}
   
   For the detailed proof of this result we refer the reader to the authors' paper \cite{gpy}, as well as to their expository account \cite{gpy-expository} and to their short article with Motohashi \cite{gpy-moto}. Our aim here is to give a very rough outline of the proof. One distinctive feature of the argument is that it `only just' works, in a way that seems rather miraculous. We will endeavour to give some sense of this. We begin with two sections of background material.
   
   \ssubsection{The Elliott-Halberstam Conjecture and level of distribution}
   
   Let $q$ be a positive integer and suppose that $a$ is prime to $q$. We write 
   \[ \psi(N; a,q):= \E_{n \leq N,n \equiv a \mdsub{q}} \Lambda(n),\]
   where $\Lambda$ is the von Mangoldt function.  For constant $q$ (and in fact for $q$ growing slowly with $N$, say $q \leq (\log N)^A$ for some fixed $A$) the prime number theorem in arithmetic progressions tells us that
   \[ \psi(N;a,q) \sim 1/\phi(q).\]
   Conditional upon the GRH, we may assert the same result up to about $q \approx N^{1/2}$. The remarkable theorem of Bombieri-Vinogradov (a proof of which the reader will find in many texts on analytic number theory, such as \cite{iwaniec-kowalski}) states that something like this is true \emph{unconditionally}, provided one is prepared to average over $q$. A weak version of the theorem is that
   \begin{equation}\label{eq1.01} \sum_{q \leq Q} \max_{(a,q) = 1} | \psi(N;a,q) - \frac{1}{\phi(q)}| \ll_{A,\epsilon} \frac{1}{(\log N)^A}\end{equation}
   for any fixed $A$ and for any $Q \leq N^{1/2 - \epsilon}$.
   
   By using \emph{sieve theory} one may show that $\psi(N;a,q) \ll N/\phi(q)$, and so the LHS of \eqref{eq1.01} is trivially bounded by 
   \[  \sum_{q \leq Q} \frac{1}{\phi(q)} \approx  \log Q.\]
   The Bombieri-Vinogradov theorem permits us to save an arbitrary power of a logarithm over this trivial bound.
   
   Even conjectures on $L$-functions (such as the GRH) appear to tell us nothing about the expression \eqref{eq1.01} when $Q \gg N^{1/2}$. Nonetheless, one may make conjectures. If $\theta \in [1/2,1)$ is a parameter then we say that the primes have \emph{level of distribution $\theta$}, or that the \emph{Elliott-Halberstam conjecture} $\mbox{EH}(\theta)$ holds, if we have the bound
   
   \begin{equation}\label{eq1.02} \sum_{q \leq Q} \max_{(a,q) = 1} | \psi(N;a,q) - \frac{1}{\phi(q)}| \ll_{A,\theta} \frac{1}{(\log N)^A}\end{equation}
   for any $Q \leq N^{\theta}$. 
   
   The \emph{full} Elliott-Halberstam conjecture \cite{EH} is that $\mbox{EH}(\theta)$ holds for all $\theta < 1$. Assuming any Elliott-Halberstam conjecture $\mbox{EH}(\theta)$ with $\theta > 1/2$, Goldston, Pintz and \yildirim can prove the remarkable result that gaps between consecutive primes are infinitely often less than some absolute constant $C(\theta)$. Assuming $\mbox{EH}(0.95971)$, they prove that
   \[ \mbox{liminf}_{n \rightarrow \infty} (p_{n+1} - p_n) \leq 16\]
   (actually they prove a slightly weaker result -- the value 0.95971 comes from unpublished computations of J. Brian Conrey).
   
   It should be stressed however that it is not expected that any conjecture $\mbox{EH}(\theta)$ for $\theta > 1/2$ will be established in the near future. There are results of Bombieri, Friedlander and Iwaniec which go a little beyond the Bombieri-Vinogradov theorem in something resembling the required manner, although experts seem to be of the opinion that these results will not help to improve the bounds on gaps between primes (cf. \cite[\S 16]{aim-notes}).
   
   \ssubsection{Selberg's weights}\label{sec1.3}
   
   This is the second section of background material.
   
   In the 1940s Selberg introduced a wonderfully simple, yet powerful, idea to analytic number theory. Write $1_P$ for the characteristic function of the primes. Then if $R$ is any parameter and if $(\lambda_d)_{d \leq R}$ is an sequence with $\lambda_1 = 1$, we have the pointwise inequality
   \[ 1_P(n) \leq \big( \sum_{\substack{d | n \\ d \leq R}} \lambda_d \big)^2\] provided that $n > R$ (the proof is obvious).
   
   This provides an enormous family of majorants for the sequence of primes. In a typical application we will be interested in something like the set of primes $p$ less than some cutoff $N$, and then $R$ will be some power $N^{\gamma}$, $\gamma < 1$. In this situation Selberg's weights majorise the primes between $N^{\gamma}$ and $N$, that is to say almost all of the primes less than $N$.
   
   What weights $\lambda_d$ should one choose? This depends on the application, but a very basic application is to the estimation of $\frac{1}{y}(\pi(x + y) - \pi(x))$, the density of primes in the interval $(x,x+y]$ (the Brun-Titchmarsh problem). In discussing this problem we will also see why it is advantageous to construct a majorant for the primes, rather than work with $1_P$ itself.
   
   For any choice of weights $\lambda_d$, then, we have
   \begin{align} \nonumber \frac{1}{y}(\pi(x+y) - \pi(x)) &\leq \E_{x+1 \leq n \leq x+y} \big( \sum_{\substack{d | n \\ d \leq R}} \lambda_d \big)^2 \\ \nonumber & =  \sum_{d \leq R} \sum_{d' \leq R} \lambda_d \lambda'_d \E_{x+1 \leq n \leq x+_y} 1_{d | n}1_{d' | n} \\ \label{display} & = \sum_{d \leq R}\sum_{d' \leq R} \frac{\lambda_d \lambda_{d'}}{[d,d']} + O\big(\frac{1}{y}\sum_{d \leq R} \sum_{d' \leq R} |\lambda_d||\lambda_{d'}|\big) .\end{align}
   Let us imagine that the weights $\lambda_d$ are chosen to be $\ll y^{\epsilon}$ in absolute value (this is always the case in practice). Then the second term here is $O(R^2 y^{2\epsilon-1})$. If $R \leq y^{1/2 - 2\epsilon}$ then this is $O(y^{- \epsilon})$ and may be thought of as an error term. This is why it is advantageous (indeed essential) to work with a majorant taken over a truncated range of divisors, and not with $1_P$ itself.
   
   The first term in \eqref{display}, 
   \[ \sum_{d \leq R}\sum_{d' \leq R} \frac{\lambda_d \lambda_{d'}}{[d,d']},\]
   is a quadratic form. It may be explicitly minimised subject to the condition $\lambda_1 = 1$, giving optimal weights $\lambda^{\mbox{\scriptsize OPT}}_d$ which are independent of $x$ and $y$, and the resultant expression may then be evaluated asymptotically. In this way one obtains the well-known bound
   \[ \frac{1}{y}(\pi(x + y) - \pi(x)) \leq (2 + \epsilon)\frac{\pi(y)}{y},\]
   valid for $y > y_0(\epsilon)$.
   
   What is the optimal choice of weights $\lambda^{\mbox{\scriptsize SEL}}_d$ for the Brun-Titchmarsh problem? The precise form will not concern us here (see, for example, \cite{nathanson}). However, it may be shown that
   \[ \lambda^{\mbox{\scriptsize SEL}}_d \approx \mu(d) \frac{\log(R/d)}{\log R}.\] (For a detailed discussion, see the appendix to \cite{green-tao-selbergsieve} and the references to work of Ramar\'e therein.)
   
 We write
 \[ \lambda^{\mbox{\scriptsize GY}}_d := \mu(d) \frac{\log(R/d)}{\log R}.\]
 These weights are very natural for two reasons: their simplicity of form, and the fact that they approximate the optimal weights for the Brun-Titchmarsh problem. There is a third reason for considering them, which comes upon recalling the formula
 
 \[ \Lambda(n) = \sum_{d | n} \mu(d) \log(n/d).\]
 We see, then, that 
 \[ \Lambda_R(n) := \sum_{\substack{d | n \\ d \leq R}} \mu(d) \log(R/d)\]
 is a kind of divisor-truncated version of $\Lambda$.
   
   We have arrived at the conclusion that the function
   \[ \frac{1}{(\log R)^2}\Lambda_R^2(n) = \big( \sum_{\substack{d | n \\ d \leq R}} \lambda^{\mbox{\scriptsize GY}}_d\big)^2\]
   might be a very useful majorant for the primes. 
     
  What might we hope to do with such a majorant? By the computation leading to \eqref{display}, we see that it is possible to find an asymptotic for 
  \begin{equation}\label{eq22} \E_{N \leq n < 2N} \Lambda_R(n)^2\end{equation}
  provided that $R \leq N^{1/2 - \epsilon}$. Later on we will wish to consider more complicated expressions involving genuine primes, such as
  \begin{equation}\label{eq23} \E_{N \leq n < 2N} \Lambda'(n +2)\Lambda_R(n)^2.\end{equation}
  Here we write $\Lambda'$ for the von Mangoldt function restricted to primes (as opposed to prime powers), thus
  \[ \Lambda'(n) := \left\{ \begin{array}{ll} \log n & \mbox{if $n$ is prime} \\ 0 & \mbox{otherwise}.\end{array}\right.\] 
  The problem of evaluating \eqref{eq22} may be thought of as a kind of approximation to the twin prime problem, though we do not know of a way to relate this expression to that problem rigourously. Expanding out, we see that \eqref{eq23} is equal to
  \[ \sum_{d \leq R} \sum_{d' \leq R} \mu(d)\mu(d')\log(R/d)\log(R/d') \E_{\substack{N \leq n < 2N\\ [d,d'] | n}} \Lambda(n + 2).\]
  Now we expect that
  \begin{equation}\label{twin} \E_{\substack{N \leq n < 2N\\ [d,d'] | n}} \Lambda'(n + 2) \approx \frac{1}{\phi([d,d'])}\end{equation} if both $d$ and $d'$ are odd.
  
  The Bombieri-Vinogradov theorem clearly offers a chance of obtaining a statement to this effect \emph{on average} over $d,d' \leq R$ if $R \leq N^{1/4 - \epsilon}$, though there is certainly still work to be done as the distribution of $[d,d']$ as $d,d'$ range over $d,d' \leq R$ is not particularly uniform. For the details (which involve moment estimates for divisor functions) see \cite[\S 9]{gpy}. 
  
Once this is done we are left with the main term
\begin{equation}\label{eq333}  \sum_{\substack{d \leq R \\ d \; \mbox{\scriptsize odd}}} \sum_{\substack{d' \leq R \\ d' \; \mbox{\scriptsize odd}}} \frac{\mu(d)\mu(d')}{\phi([d,d'])}\log(R/d)\log(R/d').\end{equation}
This term (and related expressions) may all be estimated rather accurately using the standard Dirichlet series techniques of analytic number theory, whereby the sums are expressed as integrals involving products of $\zeta$-functions.

If we had $\mbox{EH}(\theta)$, that is to say if the primes had level of distribution $\theta$, one could show that \eqref{eq23} is roughly \eqref{eq333} in the wider range $R \leq N^{\theta/2 - \epsilon}$. In particular on the full Elliott-Halberstam conjecture one could work in the range $R \leq N^{1/2`- \epsilon}$, which is essentially the same as for \eqref{eq22}.
 
 Perhaps we should make a few remarks about the form of the asymptotic for \eqref{eq333}. One may in fact show that it is
 \[ \sim  \log R \prod_{p \geq 3} \big(1 - \frac{1}{(p-1)^2}\big).\]
We will see products such as this again in \S \ref{sec3}.

\ssubsection{A strategy for gaps between primes}

We now trun to a discussion of the results of Goldston, Pintz and \yildirim themselves.
From the conceptual viewpoint it is easiest to begin by discussing the very strong \emph{conditional} results proved under the assumption of the Elliott-Halberstam conjecture. We stated in the introduction that they prove
\begin{equation}\label{eq4.1a} \mbox{liminf}_{n \rightarrow \infty} (p_{n+1} - p_n) \leq 16\end{equation}
assuming $\mbox{EH}(\theta)$ for some $\theta$ less than $1$. 

In fact, a much more general result is obtained. Let $\mathcal{H} = \{h_1,\dots,h_k\}$ be a $k$-tuple of distinct integers with $h_1 < h_2 < \dots < h_k$. A generalisation of the twin prime conjecture is that there are infinitely many $n$ such that all of $n +h_1,\dots,n + h_k$ are prime unless this is ``obviously impossible for trivial reasons'', which would be the case if there is some $p$ such that $\{h_1,\dots,h_k\}$ occupy all residue classes $\md{p}$. If this is \emph{not} the case then we say that $\mathcal{H}$ is \emph{admissible}.

Goldston, Pintz and \yildirim prove the following. 

\begin{theorem}\label{thm1} Suppose that $\mbox{\textup{EH}}(\theta)$ is known for some $\theta > 1/2$. Then there is $k_0(\theta)$ with the following property. If $k \geq k_0(\theta)$ and if $\mathcal{H} = \{h_1,\dots,h_k\}$ is an admissible $k$-tuple then for infinitely many $n$ at least \emph{two} of the numbers $n + h_1,\dots, n+ h_k$ are prime.\end{theorem}

Note in particular that 
\[ \mbox{liminf}_{n \rightarrow \infty} (p_{n+1} - p_n) \leq \min_{\substack{\{h_1,\dots,h_k\} \; \mbox{\scriptsize admissible} \\ k \geq k_0(\theta)}} (h_k - h_1).\]
It turns out that $k_0(\theta)$ can be taken to be 6 for $\theta > 0.95971$, and this leads to \eqref{eq4.1a} since the $6$-tuple $\{0,4,6,10,12,16\}$ is admissible.

Here is a very general strategy for detecting primes in admissible tuples. According to \cite{gpy-expository}, this has its origins in work of Selberg and Heath-Brown. Fix a range $[N,2N)$, suppose that $0 \leq h_1 < \dots < h_k \leq N$, and let $(\mu_n)_{N \leq n < 2N}$ be arbitrary non-negative weights which are certainly allowed to depend on the $h_i$.
We will compare
\[ Q_1 := \E_{N \leq n < 2N} \mu_n\]
with
\[ Q_2^{(i)} := \frac{1}{\log 3N}\E_{N \leq n < 2N} \Lambda'(n + h_i)\mu_n,\] for $i = 1, \dots, k$. T.~Tao remarked to me that a nice way to think of this as follows: one may renormalise so that $\sum \mu_n = 1$, and then $\rho^{(i)} := Q_2^{(i)}/Q_1$ is essentially the probability that $n + h_i$ is prime if $n$ is drawn at random from the distribution $\mu$ (one might then write the expected values in the definitions of $Q_1$ and $Q_2^{(i)}$ as integrals with respect to $\mu$).

Now if one can choose the weights $\mu$ so that $\rho^{(i)} > 1/k$, we will have upon summing over $i = 1,\dots,k$ that
\[ \E_{N \leq n < 2N} \big( \sum_{i=1}^k \frac{\Lambda'(n + h_i)}{\log 3N} - 1 \big) \mu_n > 0,\]
which means that there is some $n$ such that 
\[ \Lambda'(n + h_1) + \dots + \Lambda'(n + h_k) > \log 3N.\]
For such an $n$, at least two of $n + h_1,\dots, n+h_k$ are prime. In the probabilistic language, we have essentially used the fact that if $n + h_1,\dots,n+ h_k$ are each drawn at random from $\mu$, the expected number of primes amongst these numbers is $> 1$.

\ssubsection{Choosing good weights}

We continue the discussion of the previous section.
How should the weights $\mu_n$ be chosen to optimise the factors $\rho^{(i)}$? In retrospect, one may view most of the earlier developments on gaps between primes as attempts to find good weights $\mu_n$ in this context -- see \cite[\S 4]{gpy} for further remarks on this. 

A very good choice of weights might be
\[ \mu_n := \Lambda'(n + h_1) \dots \Lambda'(n + h_k).\]
One would indeed expect that $\rho^{(i)} \approx 1$ in this case. The only problem is that we have no idea how to prove this, one particular issue being that we cannot show that $Q_1 \neq 0$ (indeed, this is equivalent to finding $n$ such that \emph{all} of $n + h_1,\dots, n+ h_k$ are prime).

We must restrict ourselves to weights $\mu_n$ for which it is possible to estimate $Q_1$ and $Q_2^{(i)}$. As we saw in \S \ref{sec1.3}, there is a rather large class of such weights. In fact if 
\[ \mu_n = \big(\sum_{\substack{d | n \\ d \leq R}}\lambda_d\big)^2\]
then our chances are good if $R \leq N^{\theta/2 - \epsilon}$, where $\theta$ is the best exponent for which we know $\mbox{EH}(\theta)$. Crucially, there is a rather more general class of weights one is able to consider, and that is weights of the form
\begin{equation}\label{eq33b} \mu_n = \big(\sum_{\substack{d | F(n) \\ d \leq R}} \lambda_d \big)^2,\end{equation}
where $F(n) = (n+a_1)\dots (n + a_m)$ is an integer polynomial. It is an interesting exercise to reprise the arguments of \S \ref{sec1.3} in this more general context. In place of the trivial estimate
\[ \E_{n \leq N, [d,d'] | n} 1 = \frac{1}{[d,d']} + O(1)\]
which we used to derive \eqref{display} one must instead input information about
\[ \E_{n \leq N, [d,d'] | F(n)} 1.\]
But one knows that (for example) $p | F(n)$ precisely if $n \equiv -a_j \md{p}$ for some $j \in \{1,\dots,m\}$, and so such a task is not too frightening. The same is true of sums, such as \eqref{twin}, which involve the von Mangoldt function.

For the remainder of the discussion we will narrow down our search for good weights $\mu_n$ to those having the form \eqref{eq33b}. We will assume that for any sensible choice of $F$ and the $\lambda_d$ we may evaluate $Q_1$ and $Q_2^{(i)}$ for $R \leq N^{\theta/2 - \epsilon}$ using the standard techniques which we sketched in \S \ref{sec1.3}.

Now we remarked that an ideal choice of weights is 
\[ \mu_n = \Lambda'(n + h_1) \dots \Lambda'(n + h_k),\]
but we cannot compute with this choice. A closely related choice of weights is
\[ \mu_n = \Lambda_k \big((n + h_1) \dots (n + h_k)\big).\]
Here the function $\Lambda_k$ is defined by
\[ \Lambda_k(n) := \sum_{d | n} \mu(d) \log(n/d)^k,\] and so in particular $\Lambda_1 = \Lambda$. For general $k$ the function $\Lambda_k$ is supported on those integers with at most $k$ distinct prime factors. (One way to check this is to use the identity
\[ \Lambda_k = L\Lambda_{k-1} + \Lambda \ast \Lambda_{k-1},\]
where $L(n) := \log n$ and $\ast$ denotes Dirichlet convolution.)

Now in \S \ref{sec1.3} we saw the advantages of replacing $\Lambda$ with $\Lambda_R$, a divisor-truncated version of it. By analogy one might consider
\[ \Lambda_{k,R}(n) := \sum_{\substack{ d | n \\ d \leq R}} \mu(d) \log(R/d)^k.\]
This could be negative, but its square $\Lambda^2_{k,R}$ certainly cannot. Furthermore 
\[ \mu_n := \Lambda_{k,R}^2((n + h_1) \dots (n + h_k))\] is of the form \eqref{eq33b} (with $F(n) = (n + h_1)\dots (n + h_k)$). This, then, is a very natural choice of the weights $\mu_n$ and it is with some anticipation that we await the results of computing $Q_1$ and the $Q_2^{(i)}$, and hence the factors $\rho^{(i)}$. What we obtain is this:
\[ \rho^{(i)} \sim \frac{2}{k+1} \cdot \frac{\log R}{\log N}.\]
This is something of a disappointment, since it is not greater than $1/k$ even when $R = N^{\theta/2 - \epsilon}$ for $\theta$ very close to $1$ (that is, with a very strong form of the Elliott-Halberstam conjecture).

Astonishingly, a seemingly small change tips the balance in our favour. We consider instead
\[ \mu_n := \Lambda_{k+l,R}^2 ((n + h_1) \dots (n + h_k)),\]
where $0 < l \ll k$ is a further parameter. In the probabilistic language, $\rho^{(i)}$ may then be thought of (very roughly) as something like the probability that $n + h_i$ is prime given that $(n + h_1) \dots (n + h_k)$ has at most $k + l$ prime factors.

With this choice of weights it is possible to compute that
\begin{equation}\label{rho-eval} \rho^{(i)} \sim \frac{2k}{k + 2l + 1}{\frac{2l+1}{l+1}} \cdot \frac{\log R}{\log N}.\end{equation}
Note that as $k,l \rightarrow \infty$ with $l = o(k)$, this is essentially $4 \log R/k \log N$. In particular with $R := N^{\theta/2 - \epsilon}$ one has $\rho^{(i)} > 1/k$ for $k \geq k_0(\theta)$, for any $\theta > 1/2$. This is enough to establish Theorem \ref{thm1}.

\ssubsection{The unconditional result}

In this section we explain some aspects of the proof that 
\[ \mbox{liminf}_{n \rightarrow \infty} \frac{p_{n+1} - p_n}{\log n} = 0.\]
Recall that in the last section we chose weights 
\[ \mu_n := \Lambda_{k+l,R}^2((n+h_1)\dots(n + h_k)).\]
Taking $R := N^{1/4 - \epsilon}$, the quantities $Q_1$ and $Q_2^{(i)}$ can be evaluated using the Bombieri-Vinogradov theorem (rather than the Elliott-Halberstam conjecture). If $k,l$ are large with $l \ll k$ then (from \eqref{rho-eval}) the quantities $\rho^{(i)}$ are all at least $1/k - \epsilon'$ for $k \geq k_1(\epsilon)$, and hence we have that
\begin{equation}\label{almost} \E_{N \leq n < 2N} \big( \sum_{i = 1}^k \frac{\Lambda'( n + h_i)}{\log 3N}\big) \mu_n \geq (1 - \delta) \Vert \mu \Vert_1,\end{equation}
 for any $\delta > 0$, and any $k \geq k_2(\delta)$. Here
\[ \Vert \mu \Vert_1 := \E_{N \leq n < 2N} \mu_n \] is the total mass of $\mu$. This means that if $n$ is drawn at random from the distribution $\mu$, then the expected number of primes amongst the numbers $n + h_1,\dots n+h_k$ is at least $1 - \delta$. Clearly, this is not an immediately applicable result if one wishes to obtain two or more primes in the tuple $\{n + h_1,\dots, n + h_k\}$.

There is, however, a tiny bit of further information available to us. Even if $h_0 \notin \{h_1,\dots,h_k\}$, there is still a chance that if $n$ is drawn at random from $\mu$ then $n + h_0$ will be prime. One can work out that this probability is asymptotically
\[ \frac{\mathfrak{S}(h_0,h_1,\dots,h_k)}{\log N},\]
where $\mathfrak{S}(h_0,h_1,\dots,h_k)$ is a certain \emph{singular series} reflecting the arithmetic properties of the numbers $h_0,h_1,\dots,h_k$ (it is similar in form to the product defining the twin prime constant). Formally what we mean by this is that
\[ \E_{N \leq n < 2N} \frac{\Lambda'(n + h_0)}{\log 3N} \mu_n \sim \frac{\mathfrak{S}(h_0,h_1,\dots,h_k)}{\log N} \Vert \mu \Vert_1,\] and this may once again be rigorously established using the Bombieri-Vinogradov theorem.

Let $\eta > 0$ be arbitrary. Summing over all $h_0$ with 
\[ 0 \leq h_0 \leq \eta \log N,\] we obtain from \eqref{almost} that
\begin{equation}\label{eq477} \E_{N \leq n < 2N} \big( \sum_{h = 0}^{\eta \log N} \frac{\Lambda'(n + h)}{\log 3N}\big) \mu_n \geq \bigg(1 - \delta + \sum_{\substack{h_0 = 0 \\ h_0 \notin \{h_1,\dots,h_k\}}}^{\eta \log N}\frac{\mathfrak{S}(h_0,h_1,\dots,h_k)}{\log N}\bigg) \Vert \mu \Vert_1.\end{equation}
For a typical choice of $h_1,\dots,h_k$, the right-hand side will be of a predictable size. Indeed by a result of Gallagher one may infer that
\[ \sum_{\substack{h_0,\dots,h_k \leq H \\ h_i \; \mbox{\scriptsize distinct}}} \mathfrak{S}(h_0,h_1,\dots,h_k) \sim H^{k+1}\]  as $H \rightarrow \infty$. Taking expectations over all $k$-tuples $\{h_1,\dots,h_k\}$ of distinct integers with $0 \leq h_i \leq \eta \log N$, we therefore obtain from \eqref{eq477} that
\[ \E_{N \leq n \leq 2N} \big( \sum_{h = 0}^{\eta \log N} \frac{\Lambda'(n + h)}{\log 3N}\big) \mu_n \geq (1 - \delta + \eta - o(1)) \Vert \mu \Vert_1.\]
Recall that this is valid for any $\delta > 0$, provided that $k$ is sufficiently large. Taking $\delta = \eta/2$, we therefore see that (if $n$ is drawn at random from $\mu$) the expected number of primes in the interval $[n, n + \eta \log N]$ is strictly greater than $1$. In particular there is some such interval containing at least two primes.

\ssubsection{Further results}

Since the original paper of Goldston, Pintz and \yildirim several further works have appeared or are scheduled to appear giving refinements and variants of the main theorem. Here is a summary of what has been done:

\begin{itemize}

\item In the forthcoming paper \cite{gpy-precise} it is shown, by refining the ideas just described as far as seems possible, that 
\[ \liminf_{n \rightarrow \infty} \frac{p_{n+1} - p_n}{(\log p_n)^{1/2}(\log\log p_n)^2} < \infty.\]
This remarkable result asserts that the gap between primes is infinitely often almost as small as the \emph{square root} of the average gap.

\item Let $q_1 < q_2 < \dots$ be the numbers which are the product of exactly two distinct primes. Then 
\[ \liminf_{n \rightarrow \infty}(q_{n+1} - q_n) \leq 26,\] and in fact 
\[ \liminf_{n \rightarrow \infty}(q_{n+1} - q_n) \leq 6\] if one assumes the Elliott-Halberstam conjecture. These results are obtained in the paper \cite{gpyg} by the three authors and S.~W.~Graham.

\item In the original paper \cite{gpy} one also finds results concerning several primes bunching together. Thus assuming $\mbox{EH}(\theta)$ one has, for any $r \geq 2$, the bound
\[ \liminf_{n \rightarrow \infty} \frac{p_{n+r} - p_n}{\log p_n} \leq (\sqrt{r} - \sqrt{2\theta})^2.\]
One rather curious feature of this bound is that the full Elliott-Halberstam conjecture $\mbox{EH}(1)$ implies that 
\[ \liminf_{n \rightarrow \infty} \frac{p_{n+2} - p_n}{\log p_n} = 0.\]
Nothing of this sort is known for $p_{n+3} - p_n$, even conditionally.

\end{itemize}

\section{Binary digits of primes}\label{sec2}

These notes originated from a course I lectured in Part III of the Mathematical Tripos at Cambridge University in the Lent Term 2007. My aim was to work through the result of Mauduit and Rivat in the case of binary expansions of primes, and to produce a reasonably short exposition (the original paper is 49 pages long). Because we provide complete details this section is considerably more technical than either \S 1 or \S 3. Readers not interested in these technicalities may still wish to read \S \ref{sec2.3}.

\ssubsection{Statement of results}

In a recent preprint \cite{mr}, Mauduit and Rivat proved that asymptotically 50\% of the primes have odd digit sum in base 2 (and hence, of course, 50\% of the primes have \emph{even} digit sum in base 2 as well!), answering a long-standing question of Gelfond. Our aim in this section is to give a self-contained proof of their theorem in the following form, which is easily seen to imply the result as just stated.

\begin{theorem}[Mauduit-Rivat]\label{mainthm}
Let $\Lambda$ be the von-Mangoldt function, and let $s : \N \rightarrow \N$ be the function which sums the binary digits of $n$. Then 
\[ \E_{n \leq X} \Lambda(n) (-1)^{s(n)} = O(X^{-\delta})\] for some $\delta > 0$.
\end{theorem}

In fact Mauduit and Rivat proved rather more than this: they counted primes whose digit sum in base $q$ is congruent to $a \md{m}$, for any natural numbers $a,q,m$. To prove a result in this generality requires some extra technical devices and a lot more notation, so we leave the interested reader to consult the original paper \cite{mr}.

I find the result intrinsically interesting, but another reason for studying it is that it represents a pleasingly self-contained example of Vinogradov's method of `Type I and II sums' or `bilinear forms' for handling prime number sums.

In contrast with the last chapter we provide a fairly complete technical discussion.

\ssubsection{Some notation}

Throughout this section $X$ will be thought of as a large parameter. If $Y_1$ and $Y_2$ are two quantities which depend on $X$ we write $Y_1 \gtrapprox Y_2$ if $Y_1 \gg_{\epsilon} X^{-\epsilon}Y_2$ for all $\epsilon > 0$. We use the notation $Y_1 \lessapprox Y_2$ similarly. Typically we might have $Y_1 \gg Y_2\log^{-C}X$ for some $C$, but the notation is occasionally applied in even looser situations. For example we have the bound $\tau(x) \lessapprox 1$ uniformly in $x \leq X$, where $\tau(x)$ denotes the number of divisors of $x$.

If $n \in \N$ we write $\nu_2(n)$ for the $2$-exponent of $n$, the maximal $r$ such that $2^r | n$.

The letter $c$ will always denote a small, positive, absolute constant. Different instances of the letter will not necessarily denote the same constant.

\ssubsection{Vinogradov's method of Type I/II sums}\label{sec2.3}

Suppose that $f: \N \rightarrow \C$ is a function, bounded by $1$. This section concerns a method which can be often be used to show that a sum of the form
\[ \E_{n \leq X} \Lambda(n) f(n)\] is substantially smaller than $X$. The method is particularly inclined to work when $f$ is ``far'' from being multiplicative. It is clear that such a sum \emph{cannot} be small in many cases when $f$ does have some multiplicative tendancies, for example when $f(n) = \mu(n)$ or when $f(n)$ is a Dirichlet character to small modulus. 

We will develop a develop a form of the method which includes some technical refinements particularly suited to the study of the function $f(n) = (-1)^{s(n)}$. These are rather insubstantial and consist in large part of ensuring that various cutoffs are exact powers of two. In these notes we will endow the $\sim$ symbol with a rather specific meaning: if we write $\sum_{n \sim 2^{\nu}}$ then we understand that the variable $n$ is to range over the dyadic interval $[2^{\nu-1}, 2^{\nu})$.

Here is the version of Vinogradov's method that we shall require.

\begin{proposition}[Method of Type I/II sums]\label{vinogradov}
Suppose that $f : \N \rightarrow \C$ is a function with $|f(n)| \leq 1$ for all $n$. Let $\delta \in (0,1]$ be a parameter. We say that \emph{Type I sums are $\delta$-small} if 
\[ \frac{1}{X}\sum_{m \sim 2^{\mu}} |\sum_{n \in I_m} f(mn)| \leq \delta \]
whenever $2^{\mu} \leq X^{1/100}$, $I_m \subseteq [2^{\nu-1}, 2^{\nu})$ is an interval and $2^{\mu + \nu} \leq X$. We say that \emph{Type II sums are $\delta$-small} if
\[ \frac{1}{X}\big|\sum_{m \sim 2^{\mu}}\sum_{n \sim 2^{\nu}} a_m b_n f(mn)\big| \leq \delta \] uniformly for all sequences $(a_m), (b_n)$ with $|a_m|, |b_n| \leq 1$ and for all $\mu, \nu$ such that $X^{1/100} \leq 2^{\mu}, 2^{\nu} \leq X^{99/100}$.
Suppose that $f$ is a function for which Type I and Type II sums are $\delta$-small. Then
\[ \E_{n \leq X} \Lambda(n)f(n) \ll \delta^{1/2}  \log^C X.\]
\end{proposition}
\begin{remarks}
One can take $C = 11/2$. Note that both Type I and Type II sums are trivially bounded by $X$. The important feature of the proposition, then, is that a big enough gain over these trivial bounds leads to an improvement of the trivial bound on $\E_{n \leq X} \Lambda(n)f(n)$, which is $O(\log X)$.

In our statement of the result we have made a particular choice of the the ranges of $\mu, \nu$ for which estimates for Type I/II sums are required. There is considerable flexibility in the choice of these ranges. Though this matter does not concern us here, we remark that it has apparently not been completely clarified in the literature (though see \cite[\S 8]{duke-friedlander-iwaniec} for a related discussion).
\end{remarks}

When is an attempt to use the method of Type I/II sums likely to be successful in establishing a non-trivial bound for $\E_{n \leq X} \Lambda(n)f(n)$? In general, one might hope for success when $f$ does not behave ``multiplicatively''. Certainly if $f$ is multiplicative, say $f = \chi$ for some Dirichlet character $\chi$, then by choosing $a_m = \overline{f(m)}$ and $b_n = \overline{f(n)}$ we see that the Type II sums are \emph{not} always small. A similar phenomenon persists for those $f$ which are the sum of a few multiplicative functions, for example the additive character $f(n) = e(an/q)$ with $q$ a small integer. Fortunately one can estimate $\E_{n \leq X} \Lambda(n)f(n)$ for these functions by other means, namely the explicit formula and theorems on the location of zeroes of $L$-functions.

If, on the other hand, $f$ does not exhibit significant multiplicative behaviour then one may hope that the method of Type I and II sums will work. The most classical instance of this, worked out by Vinogradov in the course of proving that every large odd number is the sum of three primes, is the case $f(n) = e(\alpha n)$ where $\alpha$ is not close to a rational $a/q$. Our task here is to handle the case $f(n) = (-1)^{s(n)}$.

\ssubsection{Proof of the method of Type I/II sums.} Even though Proposition \ref{vinogradov} is not normally stated with quite the same technical refinements that we have done, the reader familiar with the basics of this theory may prefer to skip this somewhat technical section. 

Before making a start on the proof proper, we show that the assumption that Type II sums are small implies that Type I sums are small for a much larger range of $\mu$ than we have hypothesised.

\begin{lemma}[Type II implies extended Type I]\label{typeii-consequence}
Suppose that $f : \{1,\dots,N\} \rightarrow \C$ is a function with $\Vert f \Vert_{\infty} \leq 1$ for which Type I and Type II sums are $\delta$-small. Then we have the Type I estimate 
\[ \frac{1}{X}\sum_{m \sim 2^{\mu}} \big|  \sum_{n \in I_m} f(mn) \big| \ll \delta  \log X\]
in range $2^{\mu} \leq X^{99/100}$.
\end{lemma}
\proof Since Type I sums are $\delta$-small, we may assume that $2^{\mu} \geq X^{1/100}$. It clearly suffices to prove that
\begin{equation}\label{to-bound-33}\frac{1}{X} \sum_{m \sim 2^{\mu}}\sum_{n \sim 2^{\nu}} \omega_m 1_{I_m}(n) f(mn) \ll \delta \log X\end{equation} for all choices of $\omega_m$ with $|\omega_m| = 1$.
But for the cutoff $1_{I_m}(n)$, which does not factor as a product of a function of $m$ with a function of $n$, this looks like a Type II sum. To remove the cutoff, we employ a standard ``separation of variables'' trick. By the Fourier inversion formula on $\Z$ we have
\[ 1_{I_m}(n) = \int^1_0 \widehat{1}_{I_m}(\theta) e(n\theta)\, d\theta.\] Thus the left-hand side of \eqref{to-bound-33} is equal to
\[ \frac{1}{X}\int^1_0 \big( \sum_{m \sim 2^{\mu}}\sum_{n \sim 2^{\nu}} \omega_m \widehat{1}_{I_m}(\theta)e(\theta n) f(mn) \big) \, d\theta.\]
Since Type II sums are $\delta$-small and
\[ |\widehat{1}_{I_m}(\theta)| \ll \min(X,|\theta|^{-1})\]
we see that this is at most
\[ \delta  \int^1_0 \min(X, |\theta|^{-1})\, d\theta \ll \delta  \log X.\]
This concludes the proof.\endproof

\emph{Proof of Proposition \ref{vinogradov}.} The argument is essentially that of Vaughan \cite{vaughan}, but I have followed the beautiful presentation in the book of Iwaniec and Kowalski \cite{iwaniec-kowalski}. We start with the easily-verified relation
\[ \Lambda(n) = \sum_{\substack{b,c \\ bc | n}} \Lambda(b) \mu(c).\]
Let $U := X^{1/3}$ (say), and decompose this sum as
\[ \Lambda(n) = \Lambda_{\sharp \sharp}(n) + \Lambda_{\sharp \flat}(n) + \Lambda_{\flat \sharp}(n) + \Lambda_{\flat \flat}(n),\]
where $\flat$ denotes ``large'' divisors and $\sharp$ denotes ``small'' divisors, so that for example
\[ \Lambda_{\flat \sharp}(n) := \sum_{\substack{b \geq U,c < U \\ bc | n}} \Lambda(b) \mu(c).\]
We observe that
\[ \Lambda_{\sharp \flat}(n) + \Lambda_{\sharp \sharp}(n) = \sum_{b < U} \Lambda(b) \sum_{c | \frac{n}{b}}\mu(c) = 1_{n < U}\]
whilst
\begin{align*} \Lambda_{\flat \sharp}(n) + \Lambda_{\sharp \sharp}(n) &= \sum_{\substack{c < U \\ c | n}} \mu(c) \sum_{b | \frac{n}{c}} \Lambda(b) \\ &= \sum_{\substack{c < U \\ c | n}}\mu(c)\log(n/c).\end{align*}
Thus we obtain what is essentially \emph{Vaughan's identity}, 
\begin{equation}\label{vaughan-id} \Lambda(n) = -\Lambda_{\sharp \sharp}(n) + \Lambda_{\flat \flat}(n) + 1_{n < U}\Lambda(n) + \sum_{\substack{c < U \\ c | n}}\mu(c)\log(n/c). \end{equation}
Weighting by $f(n)$ and summing over $n \leq X$ we obtain
\begin{align*} & \sum_{n \leq X} \Lambda(n) f(n) \\ &= -\sum_{n \leq X}\Lambda_{\sharp\sharp}(n) f(n) + \sum_{n \leq X} \Lambda_{\flat \flat}(n)f(n) + \sum_{n < U} \Lambda(n)f(n) + \sum_{\substack{c < U \\ c | n}}\mu(c)\sum_{n \leq X} f(n) \log(n/c)\\ &= S_1 + S_2 + S_3 + S_4.\end{align*}
Our objective is to use the assumption that Type I and II sums are small in order to bound the sums $S_i$. 

\emph{Bounding $S_1$.} We have
\[ S_1 = \sum_{b,c < U} \Lambda(b)\mu(c) \sum_{n \leq X/bc} f(bcn).\] This may be rewritten as
\begin{equation}\label{eq12} \sum_{m < U^2} \omega_m \sum_{n \leq X/m} f(mn),\end{equation} where 
\[ \omega_m := \sum_{\substack{b,c < U \\ bc = m}} \Lambda(b)\mu(c).\]
By splitting into $O(\log^2 X)$ dyadic ranges for $m$ and $n$, we see that it suffices to prove that
\begin{equation}\label{to-prove-1} \sum_{m \sim 2^{\mu}} |\omega_m| \big| \sum_{n \in I_m} f(mn)\big| \ll \delta^{1/2} X\log^{C-2} X\end{equation}
whenever $I_m \subseteq [2^{\nu-1}, 2^{\nu})$ is an interval and $2^{\mu + \nu} \leq X$. This does not quite follow from Lemma \ref{typeii-consequence} since $\omega_m$ is not bounded. We do, however have $|\omega_m| \ll \tau(m)\log X$, where $\tau$ is the divisor function, and hence since $\sum_{n \leq N} \tau(n)^2 \ll N \log^3 N$ we have
\[ \sum_{m \sim 2^{\nu}} |\omega_m|^2 \ll 2^{\mu} \log^5 X.\]
This means that for any $Q > 0$ we have the estimate
\[ \sum_{\substack{m \sim 2^{\mu} \\ |\omega_m| > Q}} |\omega_m | \ll 2^{\mu}Q^{-1}\log^5 X.\]
Splitting the sum in \eqref{to-prove-1} into two parts according as $|\omega_m|$ is or is not greater that $Q$ we obtain 
\[ \sum_{m \sim 2^{\mu}} |\omega_m| \big| \sum_{n \in I_m} f(mn)\big| \ll Q^{-1}X \log^5 X + Q \sum_{m \sim 2^{\mu}}\big| \sum_{n \in I_m} f(mn)\big|.\]
By Lemma \ref{typeii-consequence} this is bounded by
\[ Q^{-1}X\log^6 X + \delta QX\log X.\]
Choosing $Q := \delta^{-1/2}\log^{5/2}X$, we obtain a bound of the required type.

\emph{Bounding $S_2$.} The sum $S_2$ may be written as 
\[ \sum_{b, c \geq U}\sum_{t \leq X/bc} \Lambda(b) \mu(c) f(bct).\]
Changing variables and rearranging, this may be written as
\begin{equation}\label{to-split} \sum_{\substack{m \geq U \\ mn \leq X}} a_m b_n f(mn),\end{equation}
where $a_m := \Lambda(m)$ and $b_n := \sum_{c \geq U: c | n} \mu(c)$. Observing that $b_n = 0$ if $n < U$, we see that the sum over $m,n$ is covered by $O(\log^2 X)$ dyadic ranges $m \sim 2^{\mu}$, $n \sim 2^{\nu}$ with $X^{0.01} \leq 2^{\mu}, 2^{\nu} \leq X^{0.99}$. We may therefore split \eqref{to-split} into $O(\log^2 X)$ sums of the form
\begin{equation}\label{to-bound-3} \sum_{m \sim 2^{\mu}} \sum_{n \sim 2^{\nu}} 1_{I(n)}(m) a_m b_n f(mn),\end{equation} where $I(n) \subseteq [2^{\mu - 1}, 2^{\mu})$ is an interval. It suffices to show that any such sum is $\ll \delta^{1/2} X\log^{C-2} X$.

These sums look rather like Type II sums, except for the presence of the cutoff $1_{I(n)}(m)$ and the fact that the sequences $a_m, b_n$ are not bounded.

To remove the cutoff we use the same device as in the proof of Lemma \ref{typeii-consequence}, writing \eqref{to-bound-3} as
\[ \int^1_0  \big( \sum_{m \sim 2^{\mu}} \sum_{n \sim 2^{\nu}} a_m e(\theta m) \widehat{1}_{I(n)}(\theta) b_n f(mn)\big) \, d\theta.\]
Since $|\widehat{1}_{I(n)}(\theta)| \ll \min(X, |\theta|^{-1})$ we see that it suffices to prove that 
\[\sum_{m \sim 2^{\mu}}\sum_{n \sim 2^{\nu}} a'_m b'_n f(mn) \ll \delta^{1/2}X\log^{C-3}X\] uniformly for all choices of $a'_m$ with $|a'_m| \leq |a_m|$ and $|b'_n| \leq |b_n|$ for all $m,n$.

This has effectively removed the cutoff that we were worried about. It remains to deal with the non-boundedness of the sequences $(a'_m)_{m \sim 2^{\mu}}, (b'_n)_{n \sim 2^{\nu}}$. The non-boundedness of $a'_m$ is rather minor, since clearly $a'_m \leq \log X$ for all $m$. Thus we may define $a''_m := a'_m/\log X$ and reduce to proving that 
\begin{equation}\label{to-prove-5} \sum_{m \sim 2^{\mu}}\sum_{n \sim 2^{\nu}} a''_m b'_n f(mn) \ll \delta^{1/2}X\log^{C-4}X\end{equation} whenever $|a''_m| \leq 1$. To deal with $b'_n$, we note that $|b'_n| \leq \tau(n)$. Thus, by the argument used in dealing with $S_1$, we have
\[ \sum_{\substack{n \sim 2^{\nu} \\  |b'_n| \geq Q}} |b'_n| \ll 2^{\nu}Q^{-1}\log^{3}X.\]
Splitting the sum in \eqref{to-prove-5} into parts according as $|b'_n| \geq Q$ or not, we conclude as before.

\emph{Bounding $S_3$.} The sum $S_3$ is trivially bounded by $U = X^{1/3}\log X$. The $\delta$-smallness of Type I sums implies (rather vacuously) that $|f(n)| \leq \delta X$. Since we are assuming that $\Vert f \Vert_{\infty} = 1$, it follows that $\delta \geq 1/X$ and hence this term may be absorbed into the bound $\delta^{1/2}X \log^C X$ that we are trying to establish.

\emph{Bounding $S_4$.} The sum may be written as
\[ \sum_{m \leq U} \sum_{n \leq X/m} \mu(m) \log n f(mn).\]
This may be split into $O(\log^2 X)$ sums of the form
\begin{equation}\label{eq46} \sum_{m \sim 2^{\mu}} \sum_{n \in I_m} \mu(m)\log n f(mn),\end{equation}
where $I_m \subseteq [2^{\nu-1}, 2^{\nu})$. This is almost a Type I sum: only the presence of the $\log n$ term prevents it being so. This term is so smooth, however, that it may be effectively removed by ``partial summation''. To do this, simply write
\[ \log n = \int^n_1 \frac{dt}{t}\] and rearrange \eqref{eq46} as 
\[ \int^{2^{\nu}}_1 \frac{dt}{t} \sum_{m \sim 2^{\nu}} \mu(m) \sum_{n \in I_{m,t}} f(mn),\]
where $I_{m,t} \subseteq [2^{\nu - 1}, 2^{\nu})$ is an interval. The smallness of Type I sums, as given in Lemma \ref{typeii-consequence}, may now be used to see that this is bounded as required.

This completes the bounding of $S_4$, and hence the proof of Proposition \ref{vinogradov}.\endproof

\ssubsection{The sum-of-digits function and related functions}

If $n$ is a positive integer and $n = \sum_{i \geq 0} n_i 2^i$ is its binary expansion, we write 
\[ s(n) := \sum_i n_i.\] This is, of course, a finite sum. For each positive integer $k$ we consider also the truncated functions
\[ s_k(n) := \sum_{i = 0}^{k-1} n_i.\]
These functions are closely related to $s(n)$, of course; on any interval $[2^k t, 2^k (t+1))$ the functions $s_k(n)$ and $s(n)$ differ by a fixed constant. The truncated functions $s_k(n)$ are periodic with period $2^k$, however, and this makes them rather amenable to study using Fourier analysis on the finite group $\Z/2^k\Z$.

In fact we will be more interested in the oscillatory functions
\[ f(n) := (-1)^{s(n)}\]
and 
\[ f_k(n) := (-1)^{s_k(n)}.\]
The functions $s_k$ and $f_k$ may, by abuse of notation, be regarded as functions on $\Z/2^k\Z$.

\begin{definition}[Finite Fourier transform] Let $k \geq 0$ be a fixed integer. Then we define
\[ \widehat{f}_k(r) := \E_{x \in \Z/2^k\Z} f_k(x) e(-rx/2^k).\]
\end{definition}

We now proceed to establish several properties of the Fourier transform $\widehat{f}_k$ which will be useful later on. We collect these here since they are all proved in a very similar way. The reader might wish to read through the proof of the first one, then skip to the next section. She might return here as each result is required.

\begin{proposition}[$L^{\infty}$ bound]\label{fourier-prop}
We have $|\widehat{f}_k(r)|  \ll 2^{-ck}$ for all $r \in \Z/2^k\Z$.
\end{proposition}
\proof We observe that 
\begin{equation}\label{fourier-expansion} \widehat{f}_k(r) = 2^{-k}(1 - e(t))(1 - e(2t)) \dots (1 - e(2^{k-1}t)),\end{equation} where $t = r/2^k$. The result now follows by observing that 
\[ |1 - e(u)||1 - e(2u)| = 4|\sin (\pi u)\sin (2\pi u)| = 8|\cos (\pi u)| (1 - \cos^2 (\pi u))\] is maximised when $\cos^2 (\pi u) = 1/3$, and attains the value $16/3\sqrt{3}$ there. Thus, grouping terms in pairs, we see that 
\[ \widehat{f}_k(r) \ll (4/3\sqrt{3})^{k/2} \ll 2^{-k/10}.\]
This concludes the proof.\endproof

\begin{proposition}[$L^1$ bound in progressions]\label{fourier-prop-2}
Suppose that $0 \leq k' \leq k$ and that $a$ is a residue class $\mdlem{2^{k'}}$. Then
\[ \sum_{\substack{r \in \Z/2^k \Z \\ r \equiv a \mdsublem{2^{k'}}}} |\widehat{f}_k(r)| \ll 2^{( \frac{1}{2} - c)(k - k')} |\widehat{f}_{k'}(a)|.\]
\end{proposition}
\begin{remark}
For the purposes of discussion suppose that $k' = a = 0$. Then this proposition states that
\[ \sum_{r \in \Z/2^k\Z} |\widehat{f}_k(r)| \ll 2^{(\frac{1}{2}-c)k}.\]
By contrast the trivial bound (arising from Parseval's identity and an application of Cauchy-Schwarz) for this quantity would be $2^{k/2}$. 

One tends to imagine that a saving over a trivial bound can be expected whenever there is some kind of ``cancellation'', such as one might expect if the function $f_k : \Z/2^k\Z \rightarrow \{-1,1\}$ were chosen at random. This setting is rather different, and this represents perhaps the most remarkable feature of the paper \cite{mr}. It is possible to show that if $f_k$ is chosen randomly then
\[ \sum_{r \in \Z/2^k\Z} |\widehat{f}_k(r)| \gg 2^{k/2}.\] Our proposition therefore captures a very specific property of the functions $f_k(n) = (-1)^{s_k(n)}$.
\end{remark}

\proof Write $S(a,k)$ for the sum on the left.
For any $r \in \Z/2^k\Z$ the expansion \eqref{fourier-expansion} yields
\[ |\widehat{f}_k(r)| = \half|1 - e(r/2^k)||\widehat{f}_{k-1}(r)|.\]
Suppose that $k \geq k' + 1$. Since $\widehat{f}_{k-1}(r)$ is periodic with period $2^{k-1}$, we may split $S(a,r)$ as a sum over two ranges $0 \leq r < 2^{k-1}$ and $2^{k-1} \leq r < 2^k$, thereby obtaining
\[ S(a,k) \leq S(a,k-1)\sup_{t \in [0,1]}\half\big(|1 + e(t)| + |1 - e(t)|\big) .\]
Now a simple exercise in Euclidean geometry confirms that
\begin{equation}\label{rat-bound} |1 + e(t)| + |1 - e(t)| \leq 2\sqrt{2},\end{equation}
with equality if and only if $t = \pm 1/4$. 
Hence
\begin{equation}\label{bound-1} S(a,k) \leq \sqrt{2} S(a,k-1).\end{equation}
This does not, in itself, suffice to establish a nontrivial bound for $S(a,k)$.

If $k \geq k' + 2$ one may improve things by splitting the sum $S(a,k)$ into \emph{four} ranges $j 2^{k-2} \leq r < (j+1)2^{k-2}$ ($j = 0,1,2,3$). This leads to 
\begin{equation}\label{bound-5} S(a,k) \leq  S(a,k-2)\sup_{t \in [0,1]} \textstyle \frac{1}{4}\displaystyle\sum_{j = 0}^3 |1 - e(t + j/4)||1 - e(2(t + j/4))|.\end{equation}
Now two applications of \eqref{rat-bound} confirm that
\begin{align*}
\sum_{j = 0}^3 & |1 - e(t + j/4)||1 - e(2(t + j/4))| \\ & = |1 - e(2t)|\big( |1 - e(t)| + |1 + e(t)|\big) + |1 + e(2t)|\big(|1 - e(t + 1/4)| + |1 + e(t + 1/4)|\big) \\ & \leq 2\sqrt{2}\big( |1 - e(2t)| + |1 + e(2t)|\big) \\ & \leq 8.
\end{align*}
Furthermore equality can never occur, and so by compactness this expression is bounded by $8 - c$ for some $c > 0$. Comparing with \eqref{bound-5} we see that
\[ S(a,k) \leq (2 - c)S(a,k-2).\] To complete the proof of the proposition one may apply this repeatedly, followed perhaps by one application of \eqref{bound-1}, to bound $S(a,k)$ above in terms of $S(a,k')$.\endproof

\ssubsection{Analysis of Type I sums}

\begin{proposition}[Estimate for Type I sums]\label{typei-prop}
Suppose that $2^{\mu} \leq X^{1/20}$. For each $m \sim 2^{\mu}$, suppose that $I_m \subseteq (2^{\nu-1}, 2^{\nu}]$ is an interval. Then 
\[ \sum_{m \sim 2^{\mu}} |\sum_{n \in I_m} f(mn)| \ll_{\epsilon} X^{19/20 + \epsilon}.\]
\end{proposition}
\begin{remark}
In \cite{mr} an alternative argument (related to the large sieve) is used, in which the same estimate is established under the weaker assumption that $2^{\mu} = O_{\epsilon}(X^{1/2 - \epsilon})$.
\end{remark}
\proof For $m,n$ in the ranges under consideration we have have $f(mn) = f_k(mn)$, where $k := \mu + \nu$. Fix a value of $m$, and write $S_m := \sum_{n \in I_m} f(mn)$. We thus have
\[ S_m = \sum_{x \in \Z/2^k \Z} f(x) 1_P(x),\]
where $P = P_m := \{mn : n \in I_m\}$.
By Parseval's identity this is 
\[ 2^k\sum_{r \in \Z/2^k \Z} \widehat{f}_k(r) \widehat{1}_P(r),\]
which is bounded by $2^k\Vert \widehat{f}_k \Vert_{\infty} \Vert \widehat{1}_P \Vert_1$. By summing a geometric series we have
\[ |\widehat{1}_P(r)| \ll \min(1, \Vert r/2^k \Vert^{-1}),\] and therefore $\Vert \widehat{1}_P \Vert_{1} \ll k \ll \log X$. It follows from Lemma \ref{fourier-prop} that \[ S_m \ll 2^k X^{-1/10}\log X.\]
Summing over $m \sim 2^{\mu}$, we obtain the required bound.\endproof

\ssubsection{Two diophantine lemmas}

In this section we give two results of a `diophantine' nature which we will use in the next section, in which Type II sums are analysed. Both of these are more-or-less standard in this subject.

\begin{lemma}[Vinogradov]\label{vinogradov-diophantine}
Suppose that $q,Q,R$ are all natural numbers less than $X$, that $\beta \in \R$, and suppose that $(a,q) = 1$. Then
\[ \sum_{x = 0}^R \min \big( Q, \Vert ax/q + \beta\Vert^{-1} \big) \lessapprox Q + q + R + QR/q.\]
\end{lemma}
\proof The key observation is that as $x$ ranges over any interval of length $q$, the numbers $ax/q \md{1}$ range over the set $0,1/q,\dots, (q-1)/q$. It follows easily that if $I$ is any interval of length at most $q$ then
\[ \sum_{x \in I} \min\big(Q, \Vert ax/q + \beta \Vert^{-1} \big) \ll Q + q\log q.\]
Since the range $x = 1,\dots,R$ may be split into at most $R/Q + 1$ such ranges, the result follows immediately.\endproof

\begin{lemma}[Equidistribution lemma]\label{eq-dist}
Suppose that $\alpha \in \R$ and that $I$ is an interval of integers with $|I| = N$. Suppose that $\delta_1,\delta_2$ satisfy $\delta_1 < \delta_2/16$, and suppose that there are at least $\delta_2 N$ elements $n \in I$ for which $\Vert \alpha n\Vert_{\R/\Z} \leq \delta_1$. Suppose that $N \geq 8/\delta_2$. Then there is some $q \leq 8/\delta_2$ such that $\Vert \alpha q \Vert_{\R/\Z} \leq 4\delta_1\delta_2^{-1}N^{-1}$.
\end{lemma}
\proof By a well-known argument of Dirichlet there is some $q \leq N$ and an $a$ coprime to $q$ such that $|\alpha - a/q| \leq 1/qN$. Set $\theta := \alpha - a/q$, let $n_0 \in I$ be fixed, and set $\beta := n_0\theta$. Then for any $n \in I$ we have, by the triangle inequality,
\[ \Vert \alpha n \Vert_{\R/\Z} \geq \Vert an/q + \beta \Vert_{\R/\Z} - \Vert \theta (n - n_0)\Vert_{\R/\Z}  \geq \Vert an/q + \beta \Vert_{\R/\Z}  - 1/q,\]
and so 
\begin{equation}\label{star-8}\Vert an/q + \beta \Vert_{\R/\Z} \leq \delta_1 + \frac{1}{q}\end{equation} for at least $\delta_2 N$ values of $n \in I$.

Now as $n$ ranges over any interval of length $q$ the numbers $an/q$ range over the set $\{0,1/q,2/q,\dots\}$. Thus in any such interval the number of $n$ satisfying \eqref{star} is no more than $\delta_1 q + 2$. Now $I$ may be divided into at most $N/q + 1$ intervals of length no more than $q$, and thus the number of $n \in I$ satisfying \eqref{star-8} is bounded by
\[ (\frac{N}{q} + 1)(\delta_1q +2).\]
On the other hand, we are assuming this quantity is at least $\delta_2 N$. Since $\delta_2 \geq 4\delta_1$, $N \geq 8/\delta_2$ and $q \leq N$ this easily implies that $q \leq 8/\delta_2$.

Now the assumption of the lemma implies, by piegonholing, that $\Vert \alpha n \Vert_{\R/\Z} \leq 2\delta_1$ for at least $\delta_2 N/2$ values of $n \in \{1,\dots,N/2\}$. We have
\[ \alpha n = \frac{an}{q} + \theta n\]
where $|\theta| \leq 1/qN$. If $n \leq N/2$, it follows that either 
\[ \Vert \alpha n \Vert_{\R/\Z} \geq 1/2q \geq \delta_2/16 > \delta_1\] or else $n$ is a multiple of $q$. But for $k \leq N/2q$ we have simply
\[ \Vert \alpha kq \Vert_{\R/\Z} = |\theta kq|.\]
Thus $|\theta k q| \leq 2\delta_1$ for at least $\delta_2 N/2$ such values of $k$, and in particular for some $k_0 \geq \delta_2 N/2$. Hence $|\theta| \leq 2\delta_1/qk_0 \leq 4\delta_1/q\delta_2 N$, which implies the result. \endproof

\ssubsection{Analysis of Type II sums}

We start with an inequality which is often used in the estimation of Type II sums, van der Corput's inequality.

\emph{Van der Corput's inequality.} Van der Corput's inequality is a kind of generalisation of the Cauchy-Schwarz inequality.

\begin{lemma}[van der Corput inequality]\label{vdc}
Let $N,H$ be positive integers and suppose that $(a_n)_{n \in \{1,\dots,N\}}$ is a sequence of complex numbers. Extend $(a_n)$ to all of $\Z$ by defining $a_n := 0$ when $n \notin \{1,\dots,N\}$. Then
\[ |\sum_{n} a_n |^2 \leq \frac{N + H}{H}\sum_{|h| \leq H}\big(1 - \frac{|h|}{H}\big)\sum_{n} a_n \overline{a_{n+h}}.\]
\end{lemma}
\proof We have
\[ \sum_n a_n = \frac{1}{H} \sum_{-H < n \leq N} \sum_{h = 0}^{H-1} a_{n+h}.\] Thus, applying the Cauchy-Schwarz inequality, we have
\begin{align*}
\big|\sum_n a_n \big|^2 &= \frac{1}{H^2} \big| \sum_{-H < n \leq N} \sum_{h = 0}^{H-1} a_{n+h}\big|^2 \\ &\leq \frac{N+H}{H^2} \sum_{-H < n \leq N} \big| \sum_{h = 0}^{H-1} a_{n+h} \big|^2 \\ & = \frac{N+H}{H^2} \sum_{-H < n \leq N} \sum_{h = 0}^{H-1} \sum_{h' = 0}^{H-1} a_{n+h} \overline{a}'_{n+h},
\end{align*}
which equals the right hand side of the claimed inequality.
This concludes the proof.\endproof

\begin{proposition}[Estimate for Type II sums]\label{typeii-prop}
Suppose that $X^{1/100} \leq 2^{\mu}, 2^{\nu} \leq X^{99/100}$. Suppose that $(a_m)_{m \sim 2^{\mu}}$ and $(b_n)_{n \sim 2^{\nu}}$ are arbitrary sequences of complex numbers with $|a_m|, |b_n| \leq 1$. Then 
\[ \sum_{m \sim 2^{\mu}}\sum_{n \sim 2^{\nu}} a_m b_n f(mn) \ll X^{1-\kappa}\]
for some absolute $\kappa > 0$.
\end{proposition}
\proof We may assume without loss of generality that $\nu \geq \mu$. Suppose for a contradiction that the result is false. We write this
\[ \sum_{m \sim 2^{\mu}} \sum_{n \sim 2^{\nu}} a_m b_n f(mn) \gtrapprox 2^{\mu + \nu}.\]

By Cauchy-Schwarz we obtain
\[ \sum_{m \sim 2^{\mu}} |\sum_{n \sim 2^{\nu}} b_n f(mn)|^2 \gtrapprox 2^{\mu + 2\nu}.\]
Let $\rho := c(\mu + \nu)$, where $c > 0$ is a constant to be chosen later, and apply the van der Corput inequality (Lemma \ref{vdc}) with $H := 2^{\rho}$. It follows that 
\[ \sum_{m \sim 2^{\mu}}\sum_{n \sim 2^{\nu}}\sum_{|h| \leq H} \big(1 - \frac{|h|}{H}\big)b_n\overline{b_{n+h}} f(m(n+h))f(mn) \gtrapprox 2^{\mu + \nu +\rho}.\]
The $h = 0$ term is negligible, and so we obtain
\begin{equation}\label{star} \sum_{1 \leq |h| \leq 2^{\rho}}\sum_{n \sim 2^{\nu}} |\sum_{m \sim 2^{\mu}} f(m(n+h))f(mn)| \gtrapprox 2^{\mu + \nu + \rho}.\end{equation}
Our aim is to obtain some cancellation in the inner sum and thereby reach a contradiction. The first step is to show that $f$ can be replaced by the truncated function $f_k$, where $k := \mu + 2\rho$ (say) is not much bigger than $\mu$. 

Now we have
\[ s(m(n + h)) - s_k(m(n+h)) = s(mn) - s_k(mn)\] if $mn, m(n+h)$ lie in the same interval $[2^kt, 2^{k}(t+1)]$. Since $mh$ is much smaller than $2^k$, this will happen most of the time. In fact for it not to hold we must have $mn \in [2^kt - 2^{\mu + \rho}, 2^kt]$ for some $t$ (that is, there is a ``carry'' on adding $mh$ to $mn$). Now for $x \leq X$ the divisor function $\tau$ satisfies $\tau(x) \lessapprox 1$, and therefore the number of ``bad'' pairs $(m,n)$ is at most
\[ \sum_{t \leq 2^{\nu - 2\rho}}\sum_{l = 2^kt - 2^{\mu + \rho}}^{2^kt} \tau(l) \lessapprox 2^{\mu + \nu - \rho}.\] The contribution of the bad pairs to \eqref{star} is therefore negligible, and we may replace that inequality by
\begin{equation}\label{star2} \sum_{1 \leq |h| \leq 2^{\rho}}\sum_{n \sim 2^{\nu}} |\sum_{m \sim 2^{\mu}} f_k(m(n+h))\overline{f_k(mn)}| \gtrapprox 2^{\mu + \nu + \rho}.\end{equation}

Let us now focus on the inner sum
\[ E_{n,h} := \sum_{m \sim 2^{\mu}}f_k(m(n+h))f_k(mn),\] which we will study using Fourier analysis on $\Z/2^k\Z$. Expanding both copies of $f_k$ using the inversion formula, we see that 
\[ E_{n,h} = \sum_{m \sim 2^{\mu}}\sum_{r,s \in \Z/2^k \Z} \widehat{f}_k(r)\widehat{f}_k(s) e\big(\frac{rm(n+h) + smn}{2^k}\big),\]
which is bounded by
\[ \sum_{r,s \in \Z/2^k \Z}|\widehat{f}_k(r)||\widehat{f}_k(s)|\min(2^{\mu} , \Vert \frac{r(n+h) + sn}{2^k}\Vert^{-1}  ).\]
From \eqref{star2} we thus have
\begin{equation}\label{star3} \sum_{r,s \in \Z/2^k \Z}|\widehat{f}_k(r)||\widehat{f}_k(s)|\sum_{n \sim 2^{\nu}}\sum_{1 \leq |h| \leq 2^{\rho}} \min(2^{\mu} , \Vert \frac{r(n+h) + sn}{2^k}\Vert^{-1}  ) \gtrapprox 2^{\mu + \nu + \rho}.\end{equation}
Define the weights
\[ \omega_{r,s} := 2^{-\mu - \nu - \rho}\sum_{n \sim 2^{\nu}}\sum_{1 \leq |h| \leq 2^{\rho}} \min(2^{\mu} , \Vert \frac{r(n+h) + sn}{2^k}\Vert^{-1}  ).\]
Then \eqref{star3} becomes
\begin{equation}\label{star4} \sum_{r,s \in \Z/2^k \Z}|\widehat{f}_k(r)||\widehat{f}_k(s)|\omega_{r,s} \gtrapprox 1.\end{equation}

To do anything with this we need to gain a greater understanding of the weights $\omega_{r,s}$.

\begin{lemma}[Size of $\omega_{r,s}$]\label{omega-size}
Suppose that $\nu_2(r+s) = t$. Then \[ \omega_{r,s} \lessapprox 2^{-\mu} + 2^{k-t - \mu - \nu} + 2^{t-k}.\]
\end{lemma}
\proof Write $(r+s)/2^k = a/2^{k-t}$, where is odd. Thus
\[ \omega_{r,s} = 2^{-\mu - \nu - \rho} \sum_{n \sim 2^{\nu}} \sum_{1 \leq |h| \leq 2^{\rho}} \min\big(2^{\mu}, \Vert \frac{a}{2^{k-t}}n + \frac{r}{2^k} h \Vert^{-1}\big).\]
Thus
\[ \omega_{r,s} \leq 2^{-\mu - \nu}\sup_{1 \leq |h| \leq 2^{\rho}} \sum_{n \sim 2^{\nu}}\min(2^{\mu} , \Vert \frac{a}{2^{k-t}}n + \frac{r}{2^{k}}h \Vert^{-1}  )\]
 By Lemma \ref{vinogradov-diophantine} we obtain
\[ \omega_{r,s} \lessapprox 2^{-\nu} + 2^{-\mu} + 2^{k-t - \mu - \nu} + 2^{t-k}.\]
Recalling that $\nu \geq \mu$, we see that this implies the claimed bound.\endproof

Let us recall \eqref{star4}. This clearly implies that there is $t$, $0 \leq t \leq k$, such that 
\begin{equation}\label{eq790} \sum_{\substack{r,s \in \Z/2^k\Z \\ \nu_2(r+s) = t}} |\widehat{f}_k(r)||\widehat{f}_k(s)|\omega_{r,s} \gtrapprox 1.\end{equation}

Now by  Lemma \ref{omega-size}, Proposition \ref{fourier-prop-2} and Parseval's inequality (in turn) we have \begin{align*} & \sum_{\substack{r,s \in \Z/2^k \Z \\ \nu_2(r+s) = t}} |\widehat{f}_k(r)||\widehat{f}_k(s)|\omega_{r,s}\\ & \lessapprox (2^{-\mu} + 2^{k-t - \mu - \nu} + 2^{t-k}) \sum_{a \in \Z/2^{t}\Z} \sum_{r \equiv a \mdsub{2^t}} |\widehat{f}_k(r)| \sum_{s \equiv -a \mdsub{2^t}} |\widehat{f}_k(s)| \\ & \lessapprox 2^{(k-t)(1-c)}(2^{-\mu} + 2^{k-t - \mu - \nu} + 2^{t-k}) \sum_{a \in \Z/2^t \Z} |\widehat{f}_t(a)||\widehat{f}_t(-a)| \\ & \leq 2^{(k-t)(1-c)}(2^{-\mu} + 2^{k-t - \mu - \nu} + 2^{t-k}) \end{align*}
for some absolute constant $c > 0$. Recalling that $k = \mu + 2\rho$ we see that the first two terms are at most $2^{4\rho - ck}$, which is not $\gtrapprox 1$ if $\rho$ is chosen sufficiently small. Thus we see that if \eqref{eq790} holds then $2^{-c(k-t)} \gtrapprox 1$, which implies that $2^t \gtrapprox 2^k$. 

Assume, then, that this is the case. Applying Parseval's identity again (or Proposition \ref{fourier-prop-2}, but this is overkill now) we obtain, for any $\eps$,

\begin{align*}
\sum_{\substack{r,s \in \Z/2^k \Z \\ \nu_2(r+s) = t \\ \omega_{r,s} \leq 2^{-\eps k}}} |\widehat{f}_k(r)||\widehat{f}_k(s)| \omega_{r,s} & \leq 2^{-\eps k} \sum_{a \in \Z/2^{t}\Z} \sum_{r \equiv a \mdsub{2^t}} |\widehat{f}_k(r)| \sum_{s \equiv -a \mdsub{2^t}} |\widehat{f}_k(s)| \\ & \leq 2^{k - t - \eps k}  \lessapprox 2^{-\eps k}.
\end{align*}

It follows that, in \eqref{eq790}, we may restrict attention to values of $r,s$ for which $\omega_{r,s} \gtrapprox 1$, thus 

\begin{equation}\label{eq791} \sum_{\substack{r,s \in \Z/2^k\Z \\  \omega_{r,s} \gtrapprox 1}} |\widehat{f}_k(r)||\widehat{f}_k(s)|\omega_{r,s} \gtrapprox 1.\end{equation}

In the next lemma we will show that there are rather few pairs $(r,s)$ with $\omega_{r,s} \gtrapprox 1$. To do this we will make use of the as yet unexploited averaging over $h$ which occurs in the definition of $\omega_{r,s}$.

\begin{lemma}[Large values of $\omega_{r,s}$]\label{omega-2}
There are $\lessapprox 2^{\rho}$ pairs $(r,s) \in \Z/2^k\Z \times \Z/2^k\Z$ such that $\omega_{r,s} \gtrapprox 1$.
\end{lemma}
\proof We know already, by Lemma \ref{omega-size}, that we must have $\nu_2(r + s) = t$, where $2^{t} \gtrapprox 2^k$. Write $(r+s)/2^k = a/2^{k-t}$, where $a$ is odd. We have
\[ \sum_{1 \leq |h| \leq 2^{\rho}} \sum_{n \sim 2^{\nu}} \min(2^{\mu}, \Vert \frac{a}{2^{k-t}} n + \frac{r}{2^k}h \Vert^{-1}) \gtrapprox 2^{\mu + \nu + \rho}.\]
Dividing the sum over $n$ into residue classes $\md{2^{k-t}}$, of which there are $\lessapprox 1$, we see that there is $\theta$ such that 
\begin{equation}\label{star-44} \sum_{1 \leq |h| \leq 2^{\rho}} \min(2^{\mu}, \Vert \theta + \frac{r}{2^k}h \Vert^{-1}) \gtrapprox 2^{\mu + \rho}.\end{equation}
Although this looks to be of the form covered by Vinogradov's lemma (Lemma \ref{vinogradov-diophantine}), the fact that $2^{\rho} \lll 2^k$ means that more information may be gleaned by appealing instead to Lemma \ref{eq-dist}. From \eqref{star-44} it follows that for $\gtrapprox 2^{\rho}$ values of $h$, $1 \leq |h| < 2^{\rho}$, we have
\[ \Vert \theta + \frac{r}{2^k}h \Vert \lessapprox 2^{-\mu}.\]
Fixing some $h_0$ with this property and considering the numbers $h - h_0$, we may assume that $\theta = 0$. Applying Lemma \ref{eq-dist}, we see that 
\[ \Vert qr/2^k\Vert_{\R/\Z} \lessapprox 2^{-\mu - \rho}\] for some $q \lessapprox 1$. Thus if $\omega_{r,s} \gtrapprox 1$ then there is some $q \lessapprox 1$ such that $|rq
 \md{2^k}| \lessapprox 2^{\rho}$ (recall that $k = \mu + 2\rho$). There are $\lessapprox 2^{\rho}$ such values of $r$, and for each of them there are just $\lessapprox 1$ values of $s$ such that $t := \nu_2(r+s)$ satisfies $2^t \gtrapprox 2^k$.\endproof
 
It is now an easy matter to show that \eqref{eq791} cannot hold which, as we have shown, implies that Type II sums are small. Indeed using Lemma \eqref{omega-2} and Proposition \ref{fourier-prop} we obtain

\[\sum_{\substack{r,s \in \Z/2^k\Z \\  \omega_{r,s} \gtrapprox 1}} |\widehat{f}_k(r)||\widehat{f}_k(s)|\omega_{r,s} \lessapprox 2^{\rho - 2ck}.\]
This contradicts \eqref{eq791} if $\rho$ is chosen sufficiently small.\endproof

\emph{Proof of Theorem \ref{mainthm}.} Almost nothing more need be said: Theorem \ref{mainthm} is an immediate consequence of Propositions \ref{vinogradov}, \ref{typei-prop} and \ref{typeii-prop}.\endproof

\section{Patterns of primes}\label{sec3}

The aim of this section is to discuss work of the author and T.~Tao on linear equations in primes. This programme is not yet completed and what has been done so far is spread across a number of papers \cite{green-tao-longprimeaps,green-tao-u3inverse,green-tao-u3mobius,green-tao-linearprimes,green-tao-nilratner,green-tao-ukmobius}. It is also discussed in several expository articles \cite{green-icm,green-longprimeaps,tao-dichotomy,tao-coates, tao-elescorial}. 

Our aim here is to do little more than try to describe what we have proved and what we hope to prove, and to furnish the reader with some idea of how the various papers fit together. In other words we hope that this section may be used as a sort of `roadmap' for readers interested in a more detailed study of the papers. There are some conspicuous absences from the current section. Our treatment of the history of the subject is very patchy (see the article \cite{kumchev-tolev} for this), and we omit any discussion of links with the ergodic theory literature (see \cite{host-survey,kra-survey} for this). References are to the versions of the papers which were available on \texttt{www.arxiv.org} on October 1st 2007. It is quite likely that published versions will have different theorem numberings.

Suppose that $\psi_1,\dots,\psi_t : \Z^d \rightarrow \Z$ are affine-linear forms. The basic questions which motivate our work are the following.

\begin{question}[Existence of prime values]\label{q1}
 Are there values of $\vec{n} \in \Z^d$ for which these forms all take prime values?  Are there infinitely many such $\vec{n}$? 
\end{question}

\begin{question}[Asymptotics]\label{q2}
How many such $\vec{n}$ are there inside the box $[-N,N]^d$?
\end{question}

In this generality, our questions contain many of the classical questions in additive prime number theory.

\begin{enumerate}
\item When $d = 1$, $t = 2$, $\psi_1(n) = n$ and $\psi_2(n) = 2m - n$ we have the \emph{Goldbach Conjecture}: is $2m$ the sum of two primes?
\item When $d = 1$, $t = 2$, $\psi_1(n) = n$ and $\psi_2(n) = n + 2$ we have the \emph{Twin Prime Conjecture}: are there infinitely pairs of primes which differ by 2?
\item When $d = 2$, $t = 3$, $\psi_1(\vec{n}) = n_1$, $\psi_2(\vec{n}) = n_2$ and $\psi_3(\vec{n}) = m - n_1 - n_2$ ($m$ odd) we have the \emph{Ternary Goldbach Conjecture}: is $m$ the sum of three primes?
\item When $d = 2$, $t = k$ and $\psi_i(\vec{n}) = n_1 + (i-1)n_2$, $i = 1,\dots,k$ we have the question of whether there exist arithmetic progressions of length $k$ consisting entirely of primes.
\end{enumerate}

Essentially nothing is known about either Question \ref{q1} or Question \ref{q2} in cases (i) and (ii). Both Questions were answered in case (iii) some seventy years ago by Vinogradov, building on earlier work of Hardy and Littlewood. Question \ref{q1} was answered in case (iv) by the author and Tao \cite{green-tao-longprimeaps}. Question \ref{q2} in that case is much harder and was answered for $k = 3$ by Chowla and van der Corput in 1939 and for $k = 4$ by the author and Tao \cite{green-tao-u3inverse,green-tao-u3mobius,green-tao-linearprimes}. Question \ref{q2} for $k \geq 5$ is one of the main goals of our current programme of research, and we shall report on what progress has been made so far. We will not give any further separate discussion of Question \ref{q1}, as this has now been exposited in many places. Particularly recommended are the articles \cite{host-survey} and \cite{kra-survey}. See also \cite{green-longprimeaps,tao-dichotomy,tao-coates,tao-elescorial}. 

There are very natural conjectural answers to Questions \ref{q1} and \ref{q2}. It is clear that congruence conditions may result in there being no, or very few, choices of $\vec{n}$ for which all of the forms $\psi_i(\vec{n})$ are prime. A trivial example is the system $\psi_1(n) = n,\psi_2(n) = n + 7$ -- consideration of this $\md{2}$ obviously implies that $\psi_1(n)$ and $\psi_2(n)$ cannot both be prime. Congruence conditions may alter our expectations in more subtle ways too. For example if $n$ is known to be prime (and is not $2$) then one feels that $n+2$ is \emph{more} likely to be prime than a random integer of the same size, for it is already known to be odd. Pulling against this, however, is the observation that if $n \neq 3$ then $n$ is congruent to either $1$ or $2 \md{3}$, and so $n+2$ is congruent to either $0$ or $1\md{3}$, but never to $2$. On this $\md{3}$ evidence one feels that $n+2$ is \emph{less} likely to be prime than a random integer of the same size. Another obvious way in which one could fail to have any prime values among the $\psi_i(\vec{n})$ is if they cannot be simultaneously positive, for example $\psi_1(\vec{n}) = n_1 - n_2$, $\psi_2(\vec{n}) = n_2 - n_3$, $\psi_3(\vec{n}) = n_3 - n_1$.

A more profound analysis of these intuitions suggests the following conjecture. In the formulation of this conjecture we use the \emph{local} von Mangoldt functions $\Lambda_{\Z/p\Z}$, defined by
\[ \Lambda_{\Z/p\Z}(x) = \left\{ \begin{array}{ll} p/(p-1) & \mbox{if $(x,p) = 1$} \\ 0 & \mbox{otherwise}.\end{array}\right.\]

\begin{dickson}
Suppose that $d,t,N \geq 1$ are integers. 
Suppose that no two of the forms $\psi_i$ are rational multiples of one another, and that no form $\psi_i$ is constant. Write $\psi_i(\vec{n}) = l_{i1}n_1 + \dots + l_{id} n_d + b_i$ and suppose that we have $|l_{ij}| \leq L$, $|b_i| \leq LN$ for some real number $L$.
Then for any convex body $K \subseteq [-N,N]^d$ we have
\[ \sum_{\vec{n} \in K \cap \Z^d} \Lambda(\psi_1(\vec{n})) \dots \Lambda(\psi_d(\vec{n})) = \beta_{\infty}\prod_p \beta_p + o_{d,t,L}(N^d),\]
where the local factors $\beta_p$ are defined by 
\[ \beta_p := \E_{\vec{x} \in (\Z/p\Z)^d} \Lambda_{\Z/p\Z}(\psi_1(\vec{x})) \dots \Lambda_{\Z/p\Z}(\psi_d(\vec{x}))\] and 
\[ \beta_{\infty} := \vol( K \cap \psi_1^{-1}(\R_{\geq 0}) \cap \dots \cap \psi_t^{-1}(\R_{\geq 0})).\] 
\end{dickson}
\remarks We have counted primes weighted using the von Mangoldt function, as this gives tidier expressions. For an unweighted version see \cite[Conjecture 1.4]{green-tao-linearprimes}. It is quite fun to play with particular cases of the conjecture. For example with $d = 1$, $t = 2$, $\psi_1(n) = n$, $\psi_2(n) = n+2$ and $K = [0,N]$ we obtain the conjecture
\[ \sum_{n \leq N} \Lambda(n) \Lambda(n+2) = 2\prod_{p \geq 3} \big( 1 - \frac{1}{(p-1)^2}  \big) N + o(N)\] for the (weighted) number of twin primes less than or equal to $N$. The numerical value of the constant here is about 1.32. Some other examples are worked out in \cite[\S 1]{green-tao-linearprimes}.

Let us return to the specific examples (i) -- (iv) mentioned above. What makes some of these questions much easier than others? The most important parameter in determining the difficulty of Questions \ref{q1} and \ref{q2} is the \emph{complexity} of the system of forms $\{\psi_1,\dots,\psi_t\}$.

The definition of complexity involves nothing more than simple linear algebra, but it is not especially illuminating at first sight. Let $s$ be a positive integer. We say that the complexity of the system $\{\psi_1,\dots,\psi_t\}$ is at most $s$ if, for any $i \in \{1,\dots,t\}$, the forms $\{\psi_1,\dots,\psi_{i-1},\psi_{i+1},\dots,\psi_t\}$ may be divided into $s+1$ classes in such a way that $\psi_i$ is not in the affine linear span of any of them. Thus the system $\{n_1,n_1 + n_2,n_1 + 2n_2, n_1 + 3n_2\}$ has complexity at most $2$ since, quite obviously, we may remove any one form and divide the remaining three forms into singleton classes whose affine span does not contain that form. However the complexity of this system is \emph{not} at most 1: if we remove the form $n_1$ then it is impossible to divide the remaining forms $\{n_1 + n_2, n_1 + 2n_2, n_1 + 3n_2\}$ into two classes such that $n_1$ is not in the affine linear span of any class. Thus the complexity of this system is exactly two. The complexity of the system $\{n_1,n_1 + n_2, n_1 + 2n_2\}$ is one, as is the complexity of the system $\{n_1, n_2, m - n_1 - n_2\}$ for any fixed $m$. The complexity of the system $\{n, n+2\}$ is apparently undefined, since if we remove the form $n$ it is impossible to partition the singleton class $\{n+2\}$ in any way such that $n$ is not contained in the affine span of a class. In such cases we say that the complexity is infinite. The system $\{n, 2m - n\}$, $m$ fixed, arising from the Goldbach Conjecture has infinite complexity.

Roughly speaking, only systems of complexity one were understood before the recent work \cite{green-tao-u3inverse, green-tao-u3mobius, green-tao-linearprimes}. Much of the theory in the complexity one case was worked out in a paper of Balog \cite{balog}, which built upon the techniques of Vinogradov. 

Complexity is the most important measure of how difficult it is to solve Questions \ref{q1} and \ref{q2}. However there are some rather trivial examples of systems of complexity greater than one for which Question \ref{q1} can be answered. For example just using the fact that there are $\gg N/\log N$ primes less than $N$ and the Cauchy-Schwarz inequality it is possible to show that there are $\gg N^4/\log^8 N$ quadruples $(n_1,n_2,n_3,n_4)$ such that all eight of the forms
\[ \{n_1, n_1 + n_2, n_1 + n_3, n_1 + n_4, n_1 + n_2 + n_3, n_1 + n_2 + n_4, n_1 + n_3 + n_4, n_1 + n_2 + n_3 + n_4\}\] are prime. This system has complexity 2.

The reason that complexity one systems have proved amenable to attack is that they can be studied using \emph{harmonic analysis}, and in particular the circle method of Hardy and Littlewood. We are now going to discuss some ideas behind this method in a rather unorthodox way. It turns out that this is the easiest (in fact so far the only) description to generalize to systems of complexity 2 and higher.
The following formulation of the principle that `harmonic analysis governs systems of complexity 1' is established in \cite{green-tao-linearprimes}.

\begin{proposition}\label{prop55} Let $N$ be a prime.
Suppose that $\{\psi_1,\dots,\psi_t\}$ is a system of affine-linear forms of complexity 1. Suppose that $f_1, \dots, f_t : \Z/N\Z \rightarrow [-1,1]$ are functions and that 
\begin{equation}\label{hyp} |\E_{\vec{n} \in (\Z/N\Z)^d} f_1(\psi_1(\vec{n})) \dots f_t(\psi_t(\vec{n}))| \geq \delta.\end{equation}
Then for each $i \in [t]$ there is some $r \in \Z/N\Z$ such that 
\begin{equation} |\E_{n \in \Z/N\Z} f_i(n) e(-rn/N)| \gg_{\delta} 1.\label{inverse-2}\end{equation}
\end{proposition}

In words, if the functions $f_1,\dots,f_t$ behave in some way unexpectedly when evaluated along the linear forms $\psi_i$, this phenomenon can be detected by evaluating the Fourier coefficients of the $f_i$.

Proposition \ref{prop55} is proved in two stages. The first step is to establish a \emph{generalized von Neumann theorem}, which is a bound of the form
\[ |\E_{\vec{n} \in (\Z/N\Z)^d} f_1(\psi_1(\vec{n})) \dots f_t(\psi_t(\vec{n}))| \leq \inf_{i \in [t]}\Vert f_i \Vert_{U^2}.\]
Here $\Vert f \Vert_{U^2}$ denotes the \emph{Gowers $U^2$-norm} of $f$ and is defined by
\[ \Vert f \Vert_{U^2}^4 := \E_{x, h_1, h_2 \in \Z/N\Z} f(x) \overline{f(x+h_1) f(x+h_2)}f(x + h_1 + h_2).\]
Results of this type are proved using nothing more than a few applications of the Cauchy-Schwarz inequality, although the notation can get quite complicated. A simple example is given in the proof of \cite[Proposition 1.9]{green-montreal}. Foundational material on the Gowers norms (including, for example, a proof that they \emph{are} norms) may be found in \cite[Chapter 11]{tao-vu} or \cite[Chapter 5]{green-tao-longprimeaps}.

This first step in the proof of Proposition \ref{prop55} leads from the hypothesis \eqref{hyp} to the conclusion that each $\Vert f_i \Vert_{U^2}$ is at least $\delta$. To obtain the conclusion \eqref{inverse-2}, then, it suffices to establish an \emph{Inverse Theorem} for the Gowers $U^2$-norm, stating that if the $U^2$-norm of $f$ is large then $f$ correlates with a linear phase. 

The proof of this result is so short we give it here. Suppose that $f : \Z/N\Z \rightarrow [-1,1]$ is a function with $\Vert f \Vert_{U^2} \geq \delta$. Define the Fourier transform $\hat{f} : \Z/N\Z \rightarrow \C$ by
\[ \hat{f}(r) := \E_{n \in \Z/N\Z} f(n) e(-rn/N).\]
Using orthogonality relations one may easily check that 
\begin{equation}\label{u2-l4} \Vert \hat{f} \Vert_4 := \big( \sum_{r\in \Z/N\Z} |\hat{f}(r)|^4 \big)^{1/4} = \Vert f \Vert_{U^2}.\end{equation}
Thus 
\[ \Vert \hat{f} \Vert_{\infty}^2 \Vert \hat{f} \Vert_2^2 \geq \Vert \hat{f} \Vert_4^4 \geq \delta^4.\]
But by Parseval's identity we have
\[ \Vert \hat{f} \Vert_2 = \Vert f \Vert_2 \leq 1,\] and so we conclude that 
\[ \Vert \hat{f} \Vert_{\infty} \geq \delta^2.\]
That is, there is some $r \in \Z/N\Z$ such that 
\begin{equation}\label{tagged-2} |\E_{n \in \Z/N\Z} f(n) e(-rn/N) | \geq \delta^2.\end{equation}
We note that \eqref{u2-l4} also gives a converse result: \emph{if} there is some $r \in \Z/N\Z$ such that 
\[ |\E_{n \in \Z/N\Z} f(n)e(-rn/N)| \geq \delta\] then $\Vert f \Vert_{U^2} \geq \delta$.

Proposition \ref{prop55} is a somewhat convincing way to formulate the idea that `harmonic analysis can handle systems of complexity one'. For a variety of reasons, however, it is not immediately applicable to an understanding of the quantity
\begin{equation}\label{eq33} \sum \Lambda(\psi_1(\vec{n})) \dots \Lambda(\psi_d(\vec{n}))\end{equation} appearing in Dickson's Conjecture.
One obvious point is that the average in Proposition \ref{prop55} is over $(\Z/N\Z)^d$, rather than over $[N]^d$ or over the lattice points inside a convex body. This is a purely technical distinction, and is dealt with in \cite[Appendix C]{green-tao-linearprimes}.

A more serious problem arises from consideration of how Proposition \ref{prop55} might be applied. There is certainly no mileage to be gained from applying it in the most na\"{\i}ve way, that is to say with $f_1 = \dots = f_t = \Lambda$. Indeed in that case condition \eqref{inverse-2} does hold (with $r = 0$), and so we cannot rule out the possibility of \eqref{hyp} holding (and, of course, Dickson's Conjecture predicts that \eqref{hyp} does hold much of the time). To eliminate the possibility of \eqref{inverse-2} holding with $r = 0$ we might split $\Lambda = 1 + (\Lambda -1)$, expand \eqref{eq33} as a sum of $2^t$ terms, and use Proposition \ref{prop55} to show that all of the terms except the one with $f_1 = \dots = f_t = 1$ are `negligible'. This also fails to work, for the simple reason that those other terms may not be negligible. Indeed if they were then the arithmetic constant in Dickson's Conjecture would be simply 1 rather than the product $\beta_{\infty}\prod_p \beta_p$ reflecting the distribution of primes $\md{2},\md{3}$ etc. 

In all of the the work of the author and Tao on primes these issues are bypassed by means of the so-called $W$-trick. The trick is easiest to describe in the context of our paper \cite{green-tao-longprimeaps} establishing the existence of arbitrarily long progressions of primes. Look at the primes $2,3,5,7,\dots$. They are very irregularly distributed $\md{2}$. However if one deletes the element $2$ and rescales the remaining primes by the map $x \mapsto (x-1)/2$ one ends up with the sequence $1,2,3,5,6,8,\dots$. This is now quite regularly distributed $\md{2}$, because there are roughly the same number of primes congruent to $1\md{4}$ as there are primes congruent to $3\md{4}$. Furthermore if one finds a long arithmetic progression in the new sequence it translates immediately to a long progression in the primes.

Unfortunately, however, this new sequence is not well-distributed $\md{3}$, as the only element congruent to $1\md{3}$ that it contains is $1$. However we could pick out the elements divisible by $3$, that is to say $3,6,9,15,\dots$ and divide through by $3$ to obtain the new sequence $1,2,3,5,\dots$. This is now well-distributed modulo both $2$ and $3$. 

The process may be continued with the primes up to some threshold $w(N)$. Choosing $w(N)$ to tend to infinity with $N$, the resulting sequences have no appreciable biases in small modulo small numbers.

It is in fact quite easy to formalise this idea and to see how it may be applied to quite general problems such as \eqref{eq33}. The sequences that result from the sieving process we have just described are of the form
\[ \{n : Wn + b \quad \mbox{is prime}\}\] where 
$W := 2 \times 3 \times \dots \times w(N)$ (hence the name `$W$-trick') and $(b,W) = 1$. For consideration of the weighted sum appearing in \eqref{eq33} it is rather natural to introduce the functions
\[ \Lambda_{b,W}(n) := \frac{\phi(W)}{W}\Lambda(Wn + b),\] which have average value roughly 1. Roughly speaking, the sum in \eqref{eq33} may be split into $\phi(W)^t$ sums of the form
\begin{equation}\label{eq34}
\sum \Lambda_{b_1,W}(\tilde{\psi}_1(\vec{n})) \dots \Lambda_{b_t,W}(\tilde{\psi}_d(\vec{n}))
\end{equation}
(for the details of this decomposition see \cite[Chapter 5]{green-tao-linearprimes}).
One might then attempt to evaluate each of these by splitting 
\[ \Lambda_{b_i,W} = 1 + (\Lambda_{b_i,W} - 1)\] and then apply Proposition \ref{prop55} to show that all of the terms except that with $f_1 = \dots = f_t = 1$ are negligible. 

Such an approach is promising, but there is one serious additional problem. Proposition \ref{prop55} only applied, as we stated it, to functions $f_i$ which are bounded by $1$. The functions $\Lambda_{b_i,W} - 1$ are certainly not bounded by one. Indeed, the harmonic analysis argument leading up to \eqref{tagged-2} relied on this boundedness in a rather essential way.

It turns out that there \emph{is} a version of Proposition \ref{prop55} which applies to functions which are not necessarily bounded by 1; this is one of the main results of \cite{green-tao-linearprimes}, and the key idea was also an important component of our earlier work \cite{green-tao-longprimeaps}.

For simplicity let us think about the von Mangoldt function $\Lambda$ itself, rather than the `W-tricked' variants $\Lambda_{b,W}$. Let $R = N^{\gamma}$, where $\gamma \in (0,1)$ is a real number, and let us recall the discussion of \S \ref{sec1}, where we observed that

\[ \frac{1}{(\log R)^2}\Lambda_R^2(n) = \big( \sum_{\substack{d | n \\ d \leq R}} \lambda^{\mbox{\scriptsize GY}}_d\big)^2\] is a sensible majorant for the characteristic function of the primes between $R$ and $N$, where 
\[ \lambda^{\mbox{\scriptsize GY}}_d := \mu(d) \frac{\log(R/d)}{\log R}.\] 
It follows the function
\[ \nu(n) := \frac{1}{\log N} \big( \sum_{\substack{d | n \\ d \leq R}} \lambda^{\mbox{\scriptsize GY}}_d\big)^2\] essentially majorizes some multiple $C_{\gamma} \Lambda(n)$ of the von Mangoldt function on the interval $[N]$ (where by \emph{essentially} we mean that there may be problems when $n \leq R$ or $n$ is a prime power, but these are highly unimportant exceptions). We note that the link we have just made to the work of Goldston, Pintz and \yildirim is by no means artificial; indeed a crucial moment in the development of \cite{green-tao-longprimeaps} occurred when Andrew Granville drew our attention to \cite{gy}.

The weights $\nu$ are much more flexible than $\Lambda$, since it is possible to evaluate such sums as 
\[ \E_{n \leq N} \nu(n) \nu(n+2)\] asymptotically by variants of the computations leading to \eqref{display}. As in those computations the key feature is that by choosing the parameter $\gamma$ to be small enough the number of terms which result from expanding out the sums over $d$ is small enough that error terms do not dominate. 

In fact (after an application of the $W$-trick discussed above and with an appropriate choice of $\gamma = \gamma(t,d)$) the weights are sufficiently flexible that they may be shown to satisfy two technical conditions called the \emph{linear forms} and \emph{correlation} conditions. These conditions were introduced in \cite[\S 3]{green-tao-longprimeaps}, and the variants of these conditions appropriate for a discussion of Dickson's Conjecture are given in \cite[\S 6]{green-tao-linearprimes}. As a result of this the weight $\nu$ qualifies to be called \emph{pseudorandom}.

As the reader may have guessed from the above discussion, it is possible to prove a version of Proposition \ref{prop55} in which the condition that the $f_i$ take values in $[-1,1]$ is relaxed to a condition $|f_i(x)| \leq \nu(x)$, where $\nu$ is a pseudorandom weight function. The first step is to establish `generalized von Neumann'-type results in which the functions $f_i$ are bounded by $\nu$, rather than just by $1$. Specifically, one deduces from an assumption \eqref{hyp} that each of the Gowers norms $\Vert f_i\Vert_{U^2}$ is somewhat large. This is once again accomplished by several applications of the Cauchy-Schwarz inequality, but the presence of the weight $\nu$ makes the details even more complicated. For a fully worked-out example, see \cite[\S 5]{green-tao-longprimeaps}. Once this is done one must establish an inverse theorem for the Gowers $U^2$-norm for functions $f$ with $|f(x)| \leq \nu(x)$. As we remarked, some new ideas are required here since the harmonic analysis argument leading up to \eqref{tagged-2} breaks down if $f$ is not bounded by one.

In fact this inverse theorem is deduced from the version with $|f_i(x)| \leq 1$ by means of the following decomposition result, which is \cite[Proposition 10.3]{green-tao-linearprimes}.

\begin{lemma}[Decomposition]\label{decomp}
Suppose that $\nu : \Z/N\Z \rightarrow \R_{\geq 0}$ is a pseudorandom measure, and that $f : \Z/N\Z \rightarrow \R$ is a function with $|f(x)| \leq \nu(x)$ for all $x$. Then we may decompose
\[ f = f_1 + f_2\]
where $|f_1(x)| \leq 1$ for all $x$ and $\Vert f_2 \Vert_{U_2} = o(1)$ as $N \rightarrow \infty$.
\end{lemma}

We shall almost nothing about the proof of this result, but it is also one of the key ingredients in \cite{green-tao-longprimeaps}.  For a discussion, see \cite[\S 6]{tao-coates}. For a broader discussion of the `energy-increment' strategy used in the proof, which appears in many different contexts in additive combinatorics, see \cite{tao-focs}.

Suppose that $|f(x)| \leq \nu(x)$ and that $\Vert f \Vert_{U^2} \geq \delta$. Applying Lemma \ref{decomp}, we see that $\Vert f_1 \Vert_{U^2} \geq \delta/2$. By the inverse theorem for the $U^2$-norm of bounded functions, there is some $r \in \Z/N\Z$ such that 
\[ |\E_{n \in \Z/N\Z} f_1(n) e(-rn/N)| \geq \delta^2/4.\]
However by the converse of the inverse theorem and the fact that $\Vert f_2 \Vert_{U^2} = o(1)$ we see that
\[ |\E_{n \in \Z/N\Z} f_2(n) e(-rn/N)| = o(1)\] (Note that the proof of this converse result did not require $f_2$ to be bounded by 1).
By the triangle inequality it follows that 
\[ |\E_{n \in \Z/N\Z} f(n) e(-rn/N)| \geq \delta^2/4 - o(1) \geq \delta^2/8,\] 
which concludes the `transference' of the inverse theorem for the $U^2$-norm from functions bounded by 1 to functions bounded by $\nu$.

Let us take stock of our position. We have indicated a proof of Proposition \ref{prop55} when the functions $f_i$ are bounded by a pseudorandom weight $\nu$, a fairly robust realisation of the principle that harmonic analysis governs the behaviour of systems of complexity one. We have split the von Mangoldt function $\Lambda$ into functions $\Lambda_{b,W}$, and rewritten the sum which interests us, namely \eqref{eq33}, as a sum of expressions of the form \eqref{eq34}. To evaluate these we effect the further splitting $\Lambda_{b,W} = 1 + (\Lambda_{b,W} - 1)$, and hope to show that any sum 
\[ \sum f_1(\tilde \psi_1(\vec{n})) \dots f_t(\tilde \psi_t(\vec{n}))\] in which at least one $f_i$ equals $\Lambda_{b,W} - 1$ is negligible. All of these functions are essentially dominated by some pseudorandom weight $\nu$ of the type considered by Goldston and Y{\i}ld{\i}r{\i}m, and so by our robust version of Proposition \ref{prop55} it suffices to rule out the possiblity that $\Lambda_{b,W} - 1$ correlates with a linear phase function; that is to say, we must establish that
\begin{equation}\label{no-linear-correlations} |\E_{n \leq N} (\Lambda_{b,W}(n) - 1) e(-rn/N)| = o(1).\end{equation}

This estimate may be established in a fairly classical fashion using the ideas of Hardy, Littlewood and Vinogradov. In \cite{green-tao-linearprimes} we proceed by first effecting some further decompositions, in keeping with our general philosophy that problems should be `transferred' to a situation where functions are bounded by 1. We skip the details (which may be found in \cite[\S 12]{green-tao-linearprimes}) and merely state that \eqref{no-linear-correlations} can be deduced from the estimate
\begin{equation}\label{mn2} |\E_{n \leq N} \mu(n)e(-rn/N)| \ll_A \log^{-A} N,
\end{equation}
for some suitably large $A$. That such an estimate holds (with \emph{any} $A$) is a well-known result of Davenport \cite{davenport}. To prove it one may use the method of Type I/Type II sums discussed in \S \ref{sec2}. In fact, Propositon \ref{vinogradov} is true with the von Mangoldt function $\Lambda$ replaced by the M\"obius function $\mu$. The proof is almost the same, relying on a decomposition of $\mu$ which is very similar to Vaughan's decomposition of $\Lambda$.

The remarks following the statement of Proposition \ref{vinogradov} are particularly relevant here. We can hope that the method of Type I/II sums will be effective in bounding \eqref{mn2} unless $r/N$ is approximately $a/q$, for some small $q$ (that is, the method ought to be successful in the `minor arc' case). Luckily in the `major arc' case one may approximate $e(-rn/N)$ by the sum of a few Dirichlet characters to modulus $q$. The resulting sums $\sum_{n \leq N} \mu(n)\chi(n)$ may then be estimating using standard techniques of analytic number theory together with information about the non-existence of zeros near $\Re s = 1$ of the $L$-functions $L(s,\chi)$: for details see \cite[Prop 5.29]{iwaniec-kowalski}.

This concludes our discussion of a proof of Dickson's Conjecture for systems of complexity one. As we have remarked, this result could also be obtained by a more classical application of the circle method. However it turns out that large parts of our discussion adapt very painlessly to systems of complexity $s > 1$, whereas this does not seem to be the case for the classical techniques.

One has, for example, the following bound of generalized von Neumann type:
\begin{equation}\label{eq47} |\E_{\vec{n} \in (\Z/N\Z)^d} f_1(\psi_1(\vec{n})) \dots f_t(\psi_t(\vec{n}))| \leq \inf_{i \in [t]}\Vert f_i \Vert_{U^{s+1}}\end{equation} for systems $\psi_1,\dots,\psi_t$ of complexity $s$,
where $\Vert f \Vert_{U^{k}}$ is the \emph{Gowers $U^k$-norm} of $f$ and is defined by
\[ \Vert f \Vert_{U^k}^{2^k} := \E_{x,h_1,\dots,h_k} \prod_{\omega_1,\dots,\omega_k \in \{0,1\}} f(x + \omega_1 h_1 + \dots + \omega_k h_k)\] (with complex conjugates being taken of the terms with an odd number of $h_i$s).
This is true even if the functions $f_i$ are only bounded by a pseudorandom weight $\nu$. A statement very close to \eqref{eq47} is proved in \cite[Appendix C]{green-tao-linearprimes}.
 
The decomposition result, Lemma \ref{decomp}, also adapts in a fairly painless manner.

The first really serious issue that we encounter is in finding a generalization of the inverse theorem for the $U^2$-norm. If $f : \Z/N\Z \rightarrow [-1,1]$ is a function such that $\Vert f \Vert_{U^3} \geq \delta$, what can be said? The most immediate difficulty with attacking this statement is the lack of a suitable formula generalizing the relation $\Vert f \Vert_{U^2} = \Vert \hat{f}\Vert_4$. A more decisive problem is revealed by consideration of the function $f(n) = e(n^2/N)$. One may check that $\Vert f \Vert_{U^3} = 1$ (this is essentially a manifestation of the fact that the third derivative of a quadratic is zero). With somewhat more effort one may check that this $f$ does not have substantial correlation with a \emph{linear} exponential $e(rn/N)$. Thus an inverse theorem for the $U^3$-norm must encode some kind of `higher harmonic analysis' which takes account of these quadratic phases as well as just linear ones. The situation is further complicated by the existence of further examples, such as $f(n) = e(n[n\sqrt{2}]/N)$, for which $\Vert f \Vert_{U^3}$ is large, but for which $f$ does not even exhibit genuine quadratic behaviour. A full discussion of examples related to these may be found in \cite{green-icm}.

It turns out that these two examples, $f(n) = e(n^2/N)$ and $f(n) = e(n[n\sqrt{2}]/N)$, can both be interpreted as objects living on a $2$-\emph{step nilmanifold}, that is to say a quotient $G/\Gamma$ where 
$G$ is a 2-step nilpotent Lie group and $\Gamma$ is a discrete and cocompact subgroup. The archetypal example is the Heisenberg example in which 
\[ G = \left(\begin{smallmatrix} 1 & \R & \R \\ 0 & 1 & \R \\ 0 & 0 & 1\end{smallmatrix}\right) \quad \mbox{and} \quad \Gamma = \left(\begin{smallmatrix} 1 & \Z & \Z \\ 0 & 1 & \Z \\ 0 & 0 & 1 \end{smallmatrix}\right).\]
Quadratic polynomials appear quite naturally in such a group $G$, for instance in the computation
\[ \left(\begin{smallmatrix} 1 & \alpha & \beta \\ 0 & 1 & \gamma \\ 0 & 0 & 1\end{smallmatrix} \right)^n = \left(\begin{smallmatrix} 1 & n\alpha & n\beta + \frac{1}{2}n(n-1)\alpha \gamma \\ 0 & 1 & n\gamma \\ 0 & 0 & 1 \end{smallmatrix} \right).\]
Taking such a sequence of matrices and postmultiplying by elements of $\Gamma$ so that all of the entries lie between $[-1/2,1/2]^3$ (that is, reducing to a fundamental domain for the action of $\Gamma$ on $G$) one soon sees the appearance of the more complicated `bracket quadratics' $\gamma n[\alpha n]$ too. A fuller discussion, with motivating remarks, may be found in \cite{green-icm}. 

What is more, correlation with an example arising in this setting is the \emph{only} way in which a function $f$ can have large $U^3$-norm. This is the inverse theorem for the $U^3$-norm, proved in \cite{green-tao-u3inverse} building on ideas of Gowers \cite{gowers-4aps}. It is conjectured that an analogous result holds for the $U^{s+1}$-norm in general, a conjecture we refer to as the \emph{Gowers Inverse conjecture} $\mbox{GI}(s)$.

\begin{conjecture}[Gowers inverse conjecture $\mbox{GI}(s)$]\label{gis}
Suppose that $f : \Z/N\Z \rightarrow [-1,1]$ is a function and that $\Vert f \Vert_{U^{s+1}} \geq \delta$.
Then there is an $s$-step nilmanifold $G/\Gamma$, a Lipschitz function $F : G/\Gamma \rightarrow [-1,1]$ and elements $g \in G$, $x \in G/\Gamma$ such that 
\[ |\E_{n \leq N} f(n) F(g^n x\Gamma)| \gg_{\delta} 1.\]
The dimension and complexity of $G/\Gamma$ and the Lipschitz constant of $F$ are all $O_{\delta}(1)$.
\end{conjecture}

The function $n \mapsto F(g^n x \Gamma)$ is called an $s$-step nilsequence. To state this conjecture properly one must of course define the notion of `complexity', and also assign a metric to $G/\Gamma$ so that the notion of Lipschitz constant may be properly formalised. A version of the conjecture was first formulated in \cite[\S 8]{green-tao-linearprimes}. There, a metric on $G/\Gamma$ was assigned quite arbitrarily. With the benefit of hindsight it is probably better to proceed as in our more recent paper \cite[\S 2]{green-tao-nilratner}, in which the notions of `complexity' and `metric' are both developed from the notion of a \emph{Mal'cev basis} for $G/\Gamma$. 

The precise statements are not important in order to understand the philosophy which lies behind Conjecture \ref{gis} and its interplay with the generalized von Neumann theorem \eqref{eq47}: it seems as though the correct `harmonics' with which to study systems $\{\psi_1,\dots,\psi_t\}$ of complexity $s$ are the $s$-step nilsequences.

Conjecture \ref{gis} is known when $s = 1$, and we proved it earlier. Note that the linear exponentials $n \mapsto e(-rn/N)$ may easily be interpreted as $1$-step nilsequences living on the nilmanifold $\R/\Z$.
As we stated, it is also known when $s = 2$, this being the main result of \cite{green-tao-u3inverse}. Tao and I hope to report progress on the general case $s \geq 3$ in the near future.

Assuming Conjecture \ref{gis}, one may start to work through the proof of Dickson's Conjecture in the complexity 1 case and attempt to generalise it to the complexity $s$ situation. Replacing occurrences of linear exponentials $e(-rn/N)$ by $s$-step nilsequences $F(g^n x\Gamma)$, the argument runs with remarkably few changes. One fairly significant extra difficulty occurs in the proof of Lemma \ref{decomp}, where a `converse' to the inverse conjecture is required. Namely, one needs to know that if 
\[ |\E_{n \leq N} f(n) F(g^n x\Gamma)| \geq \delta\] for some $s$-step nilsequence $F(g^nx\Gamma)$ then 
\[ \Vert f \Vert_{U^{s+1}} \gg_{\delta} 1,\]
where the implied constant may also depend on the complexity of $G/\Gamma$ and on the Lipschitz constant of $F$. Such a result is already present in \cite[\S 14]{green-tao-u3inverse}, and a somewhat more conceptual proof is given in \cite[Appendix E]{green-tao-linearprimes}. Both of these appendices were based on extensive conversations with members of the ergodic theory community.

By far the most serious obstacle is the last one, where we reduce to establishing a generalization of \eqref{mn2} for nilsequences. In other words we seek the bound
\[ |\E_{n \leq N} \mu(n) F(g^n x\Gamma)| \ll_A \log^{-A} N\] for all $A > 0$, where $F(g^nx\Gamma)$ is an $s$-step nilsequence arising from some $s$-step nilmanifold $G/\Gamma$, and the implied constant may also depend on the complexity of $G/\Gamma$ and the Lipschitz constant of $F$. This bound is referred to as the M\"obius and Nilsequences Conjecture $\mbox{MN}(s)$. As we remarked, the conjecture $\mbox{MN}(1)$ was essentially established by Davenport. The case $\mbox{MN}(2)$ was obtained in the paper \cite{green-tao-u3mobius}. The general case $\mbox{MN}(s)$ has recently been resolved by the authors and will appear in the short paper \cite{green-tao-ukmobius}; the key technical ingredient in that proof is the main result of \cite{green-tao-nilratner}, which may be thought of as a kind of generalization of the major-minor arc decomposition to nilsequences of arbitrary step. The method of Type I/II sums is crucial once more.

The reader wishing to study any of this work might find the following table helpful. Let me once again emphasise that the purpose of this article has been to guide potential readers through the papers \cite{green-tao-longprimeaps,green-tao-u3inverse,green-tao-u3mobius,green-tao-linearprimes,green-tao-nilratner} and \cite{green-tao-ukmobius}; we do not intend to suggest that there is no other work going on in the subject!

\cite{green-tao-longprimeaps}, \emph{The primes contain arbitrarily long arithmetic progressions,} independent of the other papers except it contains the proof of Decomposition Lemma \ref{decomp}.

\cite{green-tao-u3inverse}, \emph{An inverse theorem for the Gowers $U^3$-norm, with applications}, proof of the $\mbox{GI}(2)$ conjecture.

\cite{green-tao-u3mobius}, \emph{Quadratic Uniformity of the M\"obius function,} proof of the $\mbox{MN}(2)$ conjecture, to be largely superceded by \cite{green-tao-ukmobius} but, unlike that paper, can be understood without \cite{green-tao-nilratner}.

\cite{green-tao-linearprimes}, \emph{Linear Equations in primes,} proof that the $\mbox{GI}(s)$ and $\mbox{MN}(s)$ conjectures together imply Dickson's conjecture for systems $\{\psi_1,\dots,\psi_t\}$ of complexity $s$. The discussion in this article has been largely an exposition of some of the ideas in this paper.

\cite{green-tao-nilratner}, \emph{The Quantitative Behaviour of Polynomial Orbits on Nilmanifolds,} key technical ingredient for studying nilsequences

\cite{green-tao-ukmobius}, \emph{The M\"obius and Nilsequences Conjectures,} proof of $\mbox{MN}(s)$ conjecture for all $s$, heavily reliant on \cite{green-tao-nilratner}.

\providecommand{\bysame}{\leavevmode\hbox to3em{\hrulefill}\thinspace}
\providecommand{\MR}{\relax\ifhmode\unskip\space\fi MR }
\providecommand{\MRhref}[2]{%
  \href{http://www.ams.org/mathscinet-getitem?mr=#1}{#2}
}
\providecommand{\href}[2]{#2}

     \end{document}